\documentclass[12pt]{article}

\usepackage{amsfonts,setspace}
\usepackage{amssymb}
\usepackage{amsmath}
\usepackage{amsthm}
\usepackage{graphicx}
\usepackage[metapost]{mfpic}

\newcommand{\bt}{\begin{tabular}}
\newcommand{\et}{\end{tabular}}
\newcommand{\ba}{\begin{array}}
\newcommand{\ea}{\end{array}}

\def \Q {{\mathbb Q}}

\def \({\left(}
\def \){\right)}
\def \le {\leqslant}
\def \ge {\geqslant}
\def \proof{{\it Proof.}\ \ }
\def \ds {\displaystyle}
\def \co{{\text{\bf c}}}

\font\got = eufm10  scaled \magstep2
\font\gotbig = eufm10  scaled \magstep3
\font\gotix = eufm10  scaled \magstep0

\begin{document}

\centerline{\Large\bf Moments for generalized Farey-Brocot partitions}
\vskip+0.5cm 

\footnotetext{Supported by grants FONDECYT 7060126, RFFI 06-01-00518,
MD-3003.2006.1 (N.M.) and FONDECYT  1040975  (M.V.)}

\vskip+0.5cm
 \centerline{\bf
Nikolay Moshchevitin  \ \ and \ \  Michael Vielhaber 
 }

\vskip+0.5cm

  \vskip+1.0cm


\section{Introduction}

\subsection{Brocot partitions}

The Brocot sequences $ F_n $, $ n=0, 1, 2, \dots $  are defined as
follows.  $ F_0=\{0, 1\}=\{\frac{0}{1}, \frac{1}{1}\} $, and if
 elements of $ F_n $  are ordered in absolute value
\begin{equation}
\label{ordered_elements} 0=x_{0, n}<x_{1, n}< \dots <x_{N \left( n
  \right), n}=1, N(n)=2^{n}, 
\end{equation}
then
$$
F_{n+1} = F_n \cup Q_{n+1}, 
$$
where $ Q_{n+1} $ is defined by the rule
$$
Q_{n+1}=\{x_{i, n}\oplus x_{i-1, n}, i=1, \dots , N(n)\},  \, \, \, \, 
\frac{p}{q}\oplus \frac{p'}{q'}=\frac{p+p'}{q+q'}. 
$$
Brocot sequences were first introduced  in  \cite{Stern}, \cite{BR},
\cite{LU} (in fact, both Stern and Brocot started with
$SB_0=\left(\frac{0}{1},\frac{1}{0}\right)$, then
$SB_1=\left(\frac{0}{1},\frac{1}{1},\frac{1}{0}\right)$, so we treat
``half'' of the original sequence here).

We consider the partition of the unit interval $[0, 1]$ generated by
the points of (\ref{ordered_elements}). Let 
 $ p_{i, n}=x_{i, n}-x_{i-1, n}, i=1, \dots , N \left( n \right) $ be
the lengths of the intervals $ [x_{i-1, n}, x_{i, n}) $ . 
For  $ \beta \ge 1 $ we look for the value
$$
\sigma \left( F_n \right) =\sum_{i=1}^{N \left( n \right)}p_{i, n}^{\beta}, 
$$
the moment of order $\beta$.
The asymptotics for this sum were obtained by  Moshchevitin and
Zhigljavsky in \cite{MZ}. They proved that for $\beta >1$ the
following asymptotic formula is valid.
 \begin{equation}
\sigma_{\beta} \left( F_n \right) = \frac{2}{n^{\beta}}\frac{\zeta
  \left( 2\beta-1 \right)}{\zeta \left( 2\beta \right)} \left(1 +O
\left( \frac{\log n}{n^{ \frac{\beta - 1}{2\beta}}
 }
\right)
  \right) , 
 n \to \infty . \label{mz}
\end{equation}
Here $\zeta ( s)$ denotes the Riemann zeta function.

 The aim of the present paper is to generalize formula (\ref{mz}) to
 two special multidimensional Farey-Brocot  algorithms 
(algorithms $\hbox{\got A}$ and  $\hbox{\got B}$).

We point out that recently a better asymptotics was obtained by
Dushistova \cite{DU}. She proved that with some positive constants
$C_k= C_k (\beta), C_k^* = C_k^* (\beta )$  we have

\begin{equation}
\sigma_{\beta} \left( F_n \right) = \frac{1}{n^{\beta}} \frac{2\zeta
  \left( 2\beta-1 \right)}{\zeta \left( 2\beta \right)} +\sum_{1\le k
  < 2\beta -2}C_k\frac{1}{n^{\beta+k}} + 
\sum_{0\le k<
  \beta-2}C^*_k\frac{1}{n^{2\beta+k}} + R_\beta (n), \label{DU} 
\end{equation}
 where the remainder is
 $$
 R_\beta (n)=
 O \left( \frac{\log n}{n^{3\beta-2}}
 \right)
 $$
 in the case $2\beta \in \mathbb{Z}$, 
and
$$
 R_\beta (n) = O
\left( \frac{1}{n^{3\beta-2}}
 \right), 
 $$
when $2\beta \not\in \mathbb{Z}$. In general all summands involved
have different orders when $n\to \infty$. However,   in the case $\beta
\in\mathbb{N}$ this formula may be simplified.

\subsection{Multidimensional generalized Farey-Brocot algorithms}

There exist various multidimensional generalizations of Farey
sequences. The history starts with Hurwitz' paper \cite{HU}. As for  general
concepts and basic results we refer to  Grabiner \cite{GRA}. The
simplest  construction of multidimensional Farey nets is due to
M\"onkemeyer \cite{MON}. 
Here we introduce a general approach to  generalize the Brocot
sequence, which in some sense is similar to the construction of Farey
nets.  
In the next sections we give an analog to formula (\ref{mz})
for some special two-dimensional algorithms. 
Let ${\cal E} = \{ g_1, ..., g_{d+1}\}$ be a basis of the lattice
$\mathbb{Z}^{d+1}$.
For such a basis  ${\cal E} $  we define the cone
  $$
  {\cal C}({\cal E}) =
\{ x\in \mathbb{R}^{d+1}\, \, :\, \, \, x = \sum_{j = 1}^{d+1}
  a_jg_j, \, \, \, a_1, ..., a_{d+1} \ge 0\} .
  $$
  For a given basis
  ${\cal E} = \{g_1, ..., g_{d+1}\}$
  we can consider a natural number $K \ge 2$ and a set of integer vectors
    ${\cal E}^k_* = \{g^{k}_1, ..., g^{k}_{d+1}\}$, 
$ 1\le k \le K$
     with the following properties.
\begin{equation}
(i) \text{\ Each\ } {\cal E}^k_*  \text{\ is\ a\ basis\ for\ } 
\mathbb{Z}^{d+1}.
\text{\quad \quad \quad \quad \quad \quad \quad \quad \quad \quad
  \quad \quad \quad \quad \hfill  \quad \quad \quad \quad \quad \quad \quad
  \quad \quad \quad 
  \quad \quad \quad \quad \quad \quad \quad \quad \quad
}
\label{cond}
\end{equation}
    $(ii)$ The set of cones ${\cal C}({\cal E}^k_*), \, \, 1\le k\le K$
    forms a regular partition of the cone ${\cal C}({\cal E})$.\\ 

It means that
     ${\cal C}({\cal E})=\bigcup_{1\le k\le K}{\cal C}({\cal E}^k_*)$
    and the intersection of every  two cones ${\cal C}({\cal
    E}^{k_1}_*), {\cal C}({\cal E}^{k_2}_*)$ from this 
     union is a whole $l$-dimensional facet (for some $0\le l\le
    d$) for both  cones ${\cal C}({\cal E}^{k_1}_*)$ and 
${\cal C}({\cal E}^{k_2}_*)$.

We shall work in Euclidean space $\mathbb{R}^{d+1} $ with coordinates
$(x, y_1, ..., y_d)$.

Let the unit cube $\{z = (x, y_1, ..., y_d)\, :\, \, x=1, y_j \in [0, 1]\}$
be  partitioned into $K_0$ simplices $\Delta_k, 1\le k \le K_0,$  
in such a way that the vertices of
the simplices are among the cube's vertices. Moreover let
 the  set of vertices of each simplex
 $\Delta_k$ from this partition form a basis ${\cal E}^{0,k}$ of the lattice
$\mathbb{Z}^{d+1}$.
 A {\it generalized Farey-Brocot algorithm }  (GFBA) is a sequence
of {\it rules} $R_\nu ({\cal E}^{\nu -1, 1}, $ $\dots , {\cal E}^{\nu -1, 
  K_{\nu -1}})$ of choosing a set of bases 
$\left({\cal E}^{\nu, k} \, , \,  1\le k 
\le K_\nu({\cal E}^{\nu -1, 1}, \dots,\right.$ 
\linebreak
$\left.{\cal E}^{\nu -1, K_{\nu -1}})\right)$
for each set $({\cal E}^{\nu-1 , k}, 1\le k\le K_{\nu-1})$ 
from a previous step of the
algorithm  in such a way that every basis ${\cal E}^{\nu-1,k}$ is
decomposed into some bases from the set 
$({\cal E}^{\nu, k}, 1\le k\le K_{\nu})$ 
in such a way that 
 the conditions $(\ref{cond}) (i), (ii)$
above are satisfied.  
For given rules $R_1, ..., R_{\nu -1}$ we can construct an
infinite set of admissible rules $R_\nu$. So we can speak about a "tree" of
algorithms (compare with \cite{BALADI}). 
We shall use the gothic letter $\hbox{\got F}$
 for an individual algorithm (a precisely described set of rules).
 For the $\nu$-th set of bases we shall
use the notation ${\cal E}^{\nu , k} =\{ g^{\nu , k}_1, ..., g^{\nu , 
  k}_{d+1}\}$ and for coordinates of each vector $g^{\nu, k}_j$ we put
$g^{\nu, k}_j= (x^{\nu, k}_j, y^{\nu, k}_{j, 1}, ..., y^{\nu , 
  k}_{j, d})$.

A GFB algorithm $\hbox{\got F}$ is called {\it complete} if any
integer vector $(x, y_1, ..., y_d) \in \mathbb{Z}^{d+1}, x\ge 1, 0\le y_j
\le x, 
\, \, {\rm g.c.d.} (x, y_1, ..., y_d) = 1$ occurs as a vector from some
basis  ${\cal E}^{\nu , k}$ of the considered algorithm (there are 
GFBAs, which  are not complete).

Let  $\Theta = (1, \theta_1, ..., \theta_d), \theta _l \in [0, 1]$ be a
real vector and a GFB algorithm $\hbox{\got F}$ be given. 
To every algorithm $\hbox{\got F}$ we can construct a multidimensional
continued fraction algorithm by the following procedure. At each step
$\nu$ of the 
algorithm $\hbox{\got F}$ we choose a basis ${\cal E}^{\nu , k_\nu}
$ in such a way that that $\Theta $ can be expressed in the form $
\Theta 
= \sum_{j = 1}^{d+1}
  a_jg^{\nu , k_\nu}_j$
   and all {\it coefficients} $a_j$  of  the vector $\Theta$
  with respect to the basis ${\cal E}^{\nu, k_\nu}$ are nonnegative: $a_j \ge
  0, j = 1, ..., d+1$.
In other words, $\Theta\in {\cal C}({\cal E}^{\nu,k_\nu})$.
(One may note that in general  such a sequence of bases may not be
  unique,  for example when the 
coordinates of the vector $\Theta$ are linearly dependent over
$\mathbb{Z}$. Hence sometimes the corresponding multidimensional
continued 
fraction decomposition of the vector $\Theta$ may be not unique.)
A multidimensional continued fraction algorithm is called 
{\it  weakly convergent} in~$\Theta$~\cite{BRENTJES}, 
if for all $j, l$ 
 from the intervals 
$1\le j\le d+1, \, \, 1\le l\le d$ the sequence $y^{\nu, 
  k_\nu}_{j, l}/x^{\nu 
, k_\nu}_j$ converges to $\theta_l$. There are classical examples of
algorithms which are not weakly convergent (see 
\cite{BRENTJES}, \cite{BALADI}, \cite{SCHW}) such as Poincar\'e's  algorithm.

{\bf Lemma 1}.\, \, \, {\it Let the multidimensional continued fraction
  algorithm corresponding to the GFB algorithm $\hbox{\got F}$ be weakly
convergent. 
Then the GFBA $\hbox{\got F}$ is complete.}

(We note  that it is sufficient to suppose weak convergence of
the corresponding multidimensional continued fraction algorithm only
for vectors $\Theta$ with rational coordinates. 
Also, completeness of a GFBA does not necessarily lead to the
convergence of the  corresponding 
multidimensional continued fraction algorithm, and it is easy to
construct the corresponding example.)

\proof
Let $z = (x, y_1, ..., y_d), x\ge 1, 0\le y_j \le x, $ 
be a primitive integer point. 
Suppose it does not occur in algorithm $\hbox{\got F}$
as an element of a basis. If the described
multidimensional continued fraction algorithm is weakly convergent in the
point $\Theta =(1, y_1/x, ..., y_d/x)$ then for some sequence of basis
 ${\cal E}^{\nu , k_\nu} $ we have
$y^{\nu, k_\nu}_{j, l}/x^{\nu , k_\nu}_j \to \theta_l$ as $\nu \to\infty$
for all $j, l$. 
As we have supposed $\Theta \neq g^{\nu , k_\nu}_j$ for 
all $j$ from the interval $1\le j\le d+1$ and for all natural
$\nu$, 
this  means that for every $j$ we have $ x^{\nu , k_\nu}_j \to
+\infty $ as $\nu \to\infty$. 
But as  ${\cal E}^{\nu , k_\nu} $ is a basis and the coefficients for $\Theta$
are nonnegative  we have
 $z= x\cdot \Theta  = \sum_{j=1}^{d+1}  a_jg^{\nu , k_\nu}_j$ with
nonnegative integers  $a_j$ and at least one  of them (say $a_m$)  is $\ge 1$.
Then $x \ge a_m\cdot x^{\nu , k_\nu}_m\to +\infty$ and this is a
contradiction.~\hfill$\Box$ 

Each GFB algorithm $\hbox{\got F}$  generates a sequence of partitions 
(``tilings'') ${\rm Til }_\nu (\hbox{\got F})$ of the unit cube
$[0, 1]^d$  and a sequence of graphs $T_\nu (\hbox{\got F})$ as follows:

Let a GFB algorithm $\hbox{\got F}$ be given. We look for the set of
all bases ${\cal E}$ at the $\nu$-th step of our algorithm. The number
of such bases may vary according to {\got F}, 
but the corresponding cones ${\cal C}({\cal
  E})$ form a regular partition of the cubic cone
$\{z=(x, y_1, ..., y_d)\, :\, \, \, 
x\ge 0, y_j \in [0, 1]\}$. We restrict this partition on the set
$\{z=(x, y_1, ..., y_d)\, :\, \, \, x = 1, y_j \in [0, 1]\}$ and obtain the
partition 
${\rm Til}_\nu  (\hbox{\got F}) $ of the unit cube $[0, 1]^d$ into
simplices $\Delta={\cal C}({\cal E})\cap \{z=(x, y_1, ..., y_d)\, :\, \, \, 
x = 1, y_j \in [0, 1]\}$. 
The main object of the paper is the sum
$$
\sigma_{n, \beta } (\hbox{\got F}) = \sum_{\Delta\in {\rm Til}_n
(\mbox{\gotix F})} ({\rm mes} \Delta )^\beta,
$$
the moment of order $\beta$.

Obviously for any GFBA and for any natural $n $ we have
\begin{equation}
\sigma_{n, 1 } (\hbox{\got F}) =1. \label{one}
\end{equation}

The following simple statement is well-known (see, for example
\cite{WOO}). 

{\bf Lemma 2.}\, \, \, {\it

Let the simplex $\Delta$ correspond to the basis ${\cal E} = \{ g_1, ..., 
g_{d+1}\}$, and the vector $g_j$ from this basis have coordinates 
$g_j = (x_j, y_{j,1}, ..., y_{j,d})$. 
Then
$$ {\rm mes }\Delta =\frac{1}{d!x_1\cdots x_{d+1}}.
$$}

\proof
See \cite[Thm.~9]{WOO}. \hfill$\Box$

The graph $T_\nu (\hbox{\got F})$ is defined as follows. 
The set
$V_\nu (\hbox{\got F})$ of its vertices is the set of all vectors from
all  bases of the $\nu$-th step of the algorithm, and we have an edge between
vertices $u$ and $v$ if the integer vectors $u, v$ belong to the same
basis ${\cal E}$. 
We also consider the graph $T(\hbox{\got F})$ whose
vertices $V (\hbox{\got F})$ are the vectors of {\it  all} bases 
${\cal  E}$  appearing in the algorithm $\hbox{\got F}$ and there
exists an edge between vertices $u$ and $v$ if and only if vectors
$u, v$ belong to the same 
  basis ${\cal E}$.
We emphasize that if {\got F} is complete, 
$\left\{\left(\frac{y_1}{x},\dots,\frac{y_d}{x}\right)\ | \
(x,y_1,\dots,y_d)\in V(\mbox{\got F})\right\} =
\Q^d\cap [0,1]^d$. 
 Clearly $T_n(\hbox{\got F})$
is a subgraph of $T(\hbox{\got F})$.

 We define a GFB algorithm to be {\it finite}  if there is a
positive constant $ M(\hbox{\got F})$ such that for any vertex $v\in V
(\hbox{\got F})$ of the graph $T(\hbox{\got F})$ its degree $\deg(v)$ (the
number of edges with the endpoint in this vertex) is bounded by  $
M(\hbox{\got F})$.

{\bf Lemma 3.}\, \, \, {\it If the  GFBA $\hbox{\got F}$ is finite then it
  is complete.} 

\proof
Let $z = (x, y_1, ..., y_d), x\ge 1, 0\le y_j \le x, $ 
be a primitive integer point and 
 $\Theta =(1, y_1/x, ..., y_d/x)$. Suppose $z$ does not occur in algorithm
$\hbox{\got F}$ as an element of a basis, $z\not\in V(\mbox{\got F})$. 
As in Lemma 1,  we define for every $\nu$ a $k_\nu$ such that 
$z\in {\cal C}({\cal E}^{\nu , k_\nu })$, and look for the basis ${\cal
  E}^{\nu , k_\nu} =\{ g^{\nu , k_\nu}_1, ..., g^{\nu , k_\nu}_{d+1}\}$ 
and for  the corresponding simplex $\Delta_{\nu}$ from partition ${\rm
  Til}_\nu$.  
We have $\Theta \in \Delta_\nu, \, \, \forall \nu$.
 We shall prove that the
first coordinates $x^{\nu , k_\nu }_{j}$ 
of the basis vectors 
$g^{\nu, k_\nu }_{j}$  
tend to $+\infty$ as $\nu \to\infty$. 
Then the lemma will be proved as $ z = \sum_{j=1}^{d+1}
\lambda_{j}g^{\nu , k_\nu 
}_{j}$ with nonnegative $ \lambda_{j}$ (and one of the  $
\lambda_{j}$  must be positive). 
So $ x \ge \min_{1\le j\le d+1} x^{\nu , k_\nu }_{j}$,
 and this is a contradiction. 

To do it we must use the finiteness property of our algorithm.
Let $a_1, a_2, ..., a_{d+1}$ be vertices of the  simplex
$\Delta_{\nu}$. Now we fix a vertex  of this simplex (say $a_{d+1}$)
and show that for large 
enough $\nu '$ this vertex will not be a vertex of $\Delta_{\nu
  '}$. 
By the finiteness of the algorithm, we may assume that for $\nu ' \ge
\nu $  no 
additional edge with vertex $a_{d+1}$ appears inside simplex
$\Delta_{\nu '}$ during our algorithm. 
Also, during the algorithm only finitely many 
new edges may appear in the vertices $a_1, .., a_d$. 
Since at each
step of our algorithm we must choose a partition of each simplex into 
smaller ones, if $a_{d+1}$ is still a vertex of a
simplex $\Delta '$ from a partition ${\rm Til}_{\nu '}$ with large
$\nu '$, 
then the other vertices $a_1'= a_1' (\nu '), ..., a_d' =a_d ' (\nu ')$ of 
$\Delta '$ must lie on edges $[a_{d+1}, a_j], j = 1, .., d$  and do
not coincide with the endpoints $ a_1, ..., a_d$.  

Let ${\cal E}^{\nu ', k_{\nu '}}
=\{ g^{\nu ' , k_{\nu '}}_1, ..., g^{\nu ' , k_{\nu '}}_{d+1}\}$ be the  
corresponding basis. Then $ g^{\nu ' , k_{\nu '}}_{d+1} =g^{\nu  , 
  k_{\nu }}_{d+1} $ and for every $j\in [1, .., d]$ we have 
$g^{\nu ',k_{\nu'}}_{j} = 
\mu_j g^{\nu, k_{\nu }}_{d+1} +  g^{\nu  , k_{\nu }}_{j}$
with {\it positive } $\mu_j$. 
It means that $a_j ' (\nu ') \to a_j, 1\le j \le d$ when $\nu ' \to
\infty $. 
So for large $\nu ' $ the point
$\Theta$ will not lie in the simplex with vertex $a_{d+1} $ from the
partition 
${\rm Til }_{\nu '}$. We see that for large $\nu '$ all vertices of
the simplex $\Delta_{\nu '}$ will differ from the vertices of the
simplex  $\Delta_{\nu}$. 
Now $\ds\min_{1\le j\le d+1} g^{\nu', k_{\nu '}}_{j, 1}
> \min_{1\le l\le d+1} g^{\nu  , k_{\nu }}_{j, 1}$, and the first
coordinates  of the basis vectors $g^{\nu , k_\nu }_{j, 1}$ tend to
$+\infty$ as $\nu \to\infty$. The proof is complete. \hfill$\Box$.

For a finite GFB algorithm $\hbox{\got F}$ we consider the  Dirichlet series
$$
L(\hbox{\got F}, \beta ) = \sum_{v\in V (\hbox{\gotix F})} 
\frac{\deg(v)}{(x(v))^\beta},
$$
where $x(v)$ is the first coordinate of the integer vector
$v=(x,y_1,\dots,y_d)$. 
Since for the fixed value of $x$  the number of integer vectors of the form
$(x, y_1, ..., y_d), \, \, 0\le y_l\le x, \, \, {\rm g.c.d.} (x, y_1,
..., y_d) = 1$  is bounded by $(x+1)^d-2^d$, the series for $L(\hbox{\got
  F}, \beta )$ converges when $\beta > d+1$.

Let $ a= \left(\frac{a_1}{q}, ..., \frac{a_d}{q}\right) \in
\(\mathbb{Q}^+\)^d$  and ${\rm g.c.d.}(q, a_1, ..., a_d) =1$. For
every $a \in \(\mathbb{Q}^+\)^d$ we 
define $ q(a)=q$. 
Recall that if the GFBA is finite, then it is complete, and
hence the corresponding vertex $a=(q,a_1,\dots,a_d)$ 
occurs as a vertex of our graph $T(\hbox{\got  F})$. 
Let $\deg(a)$ be its degree in this graph. Then for $ \beta > d+1$ we have
\begin{equation}
L(\hbox{\got F}, \beta ) = \sum_{a \in \(\mathbb{Q}^+\)^d} \frac{
 \deg(a)}{(q(a))^\beta} = \sum_{q=1}^{+\infty} 
 \frac{ \displaystyle{\sum_{l=1}^{M(\hbox{\gotix
F})} l\times G_l(q)}}{q^\beta}, \label{L}
\end{equation}
 where $G_l(q)$ denotes the number of rational points  $a$ with common
 denominator $q$ such that $\deg(a)$~=~$l$.

In the case $d=1$ for the algorithm of taking medians of neighbouring
on each step (as it was described in Section 1.1.), two vertices ($0/1$ 
and $1/1$) of the graph $T$ have degree 1 and all other vertices have
degree 2. So for the classical one-dimensional algorithm the
Dirichlet series is
$$
L(\beta ) = 2\times\frac{\zeta \left( \beta-1 \right)}{\zeta \left(
  \beta \right)} = 2\times \sum_{q=1}^\infty \frac{\varphi (q)}{
  q^{\beta}} 
$$
(where $\varphi (\cdot )$ is Euler's totient function), and formulas
(\ref{mz}), (\ref{DU}) give  good asymptotics for the moments of the
partition generated by the algorithm under consideration with the
coefficient $L( 2\beta ) $ in the main term.

The main result of this paper is obtaining nice asymptotic
formulas for $\sigma_{n, \beta } (\hbox{\got F}), \, \, n \to \infty$
when $ d = 2$ and $\hbox{\got F}$ is one of two simplest finite GFBA 
algorithms. 
In Section 2 we consider algorithm $\hbox{\got A}$. It seems to be new
but it is related to the (generalized) Poincar\'e algorithm (see 
\cite{PUA}, \cite[ch.~21]{SCHW}, \cite{NOUG}). 
Also algorithm  $\hbox{\got A}$ is related to the construction from
\cite{WOO}. 
In Section 3 we consider algorithm $\hbox{\got B}$, 
which was  introduced by M\"onkemeyer in \cite{MON}. 


We note that 
both algorithms use a fixed rule for choosing the 
partition of Til$_n$ from Til$_{n-1}$, for all partitions Til$_n$.
The description of algorithm 
$\hbox{\got A}$ seems to be somewhat more complicated than that of
algorithm $\hbox{\got B}$, however the rule for the new partitions in
algorithm  
$\hbox{\got A}$  does not depend on the order of the vectors in the
preceeding  basis, while the rule for algorithm $\hbox{\got B}$ does.

\section{Algorithm $\hbox{\gotbig A}$}

\subsection{The description of Algorithm $\hbox{\got A}$}

We fix the initial partition of the unit square $\{z=(x, y_1, y_2)\, \colon \, 
x=1, y_{1, 2} \in [0, 1]\}$ into two triangles
 $\Delta^{0,1} $ with vertices $(1, 0, 0), (1, 1, 0), (1, 0, 1)  $
and $\Delta^{0,2}$ with vertices $ (1, 1, 0), (1, 0, 1), (1, 1, 1)$. Vertices
of both triangles form bases ${\cal E}^{0,1}$ and ${\cal E}^{0,2}$ 
of the integer lattice $\mathbb{Z}^3$. Now we suppose that a basis
$$
{\cal E}^{\nu,j} = \{ g^{\nu , j}_{1}, g^{\nu , j}_2, g^{\nu , j}_3\}
$$
(and the corresponding triangle $ {\Delta}^{\nu,j} $ of the partition of
the unit square $[0, 1]^2$ into triangles) 
 occurs in our algorithm, and  we define
the rule for constructing the bases for the next step of algorithm.
 In our algorithm $\hbox{\got B}$ the rule also  will be the same for
each step of the algorithm and for each basis. 
Namely, for the basis ${\cal E}^{\nu,j}$ which occurs at the $\nu$-th
step, we take $6$ bases 
${\cal E}^{\nu +1,6(j-1)+i}, \, i = 1, 2, 3, 4, 5, 6$ by the following
formulas 
$$ {\cal E}^{\nu+1,6(j-1)+1} = \{ g^{\nu , j}_{1}, g^{\nu , j}_1+g^{\nu
  , j}_2, g^{\nu , j}_1+g^{\nu , j}_3\}, 
$$
$$ {\cal E}^{\nu+1,6(j-1)+2} = \{
g^{\nu , j}_{2}, g^{\nu , j}_2+g^{\nu , j}_1, g^{\nu , j}_2+g^{\nu , j}_3\}, 
$$
$$ {\cal E}^{\nu+1,6(j-1)+3} = \{
g^{\nu , j}_{3}, g^{\nu , j}_3+g^{\nu , j}_1, g^{\nu , j}_3+g^{\nu , j}_2\}, 
$$
$$ {\cal E}^{\nu+1,6(j-1)+4} = \{
g^{\nu , j}_1+g^{\nu , j}_{2}, g^{\nu , j}_1+g^{\nu , j}_3, g^{\nu
  , j}_1+g^{\nu , j}_2+g^{\nu , j}_3\}, 
$$
$$ {\cal E}^{\nu+1,6(j-1)+5} = \{
g^{\nu , j}_2+g^{\nu , j}_{1}, g^{\nu , j}_2+g^{\nu , j}_3, g^{\nu
  , j}_1+g^{\nu , j}_2+g^{\nu , j}_3\}, 
$$
$$ {\cal E}^{\nu+1,6(j-1)+6} = \{
g^{\nu , j}_3+g^{\nu , j}_{1}, g^{\nu , j}_3+g^{\nu , j}_2, g^{\nu
  , j}_1+g^{\nu , j}_2+g^{\nu , j}_3\}. 
$$

We point out  that the construction of the set of bases ${\cal E}^{\nu
  +1,6(j-1)+i}, \, 1\le i\le 6$ does not depend on the {\it order} of
vectors in the basis ${\cal E}^{\nu ,j}$.

It is easy to see that the described rule satisfies the conditions
(\ref{cond})$(i), (ii)$ --- 
each new set of vectors  ${\cal E}^{\nu +1,6(j-1)+i}$
forms  a basis of the integer lattice and the cones 
${\cal C}({\cal E}^{\nu+1,6(j-1)+i}), \, 1\le i\le 6$
form a regular partition of the cone
${\cal C}({\cal E}^{\nu,j})$.
Obviously this algorithm is finite, and
from Lemma 3 it follows that algorithm $\hbox{\got A}$ is complete.
Hence, for any $\xi = (p, a_1, a_2)\in \mathbb{Z}^3, \, \, 
p\ge 1, 0\le a_1, a_2 \le p,  \text{g.c.d.}(p,a_1,a_2)=1$ 
there exist $m\le p$ and $j$ such that $\xi \in {\cal E}^{m,j}$.

We shall show in 2.4 that the multidimensional continued fraction algorithm 
corresponding to algorithm $\hbox{\got A}$  weakly
converges everywhere.

\subsection{Algorithm $\hbox{\got A}$ in terms of constructing rational
 points in the square $[0, 1]^2$}

For two rational points
 $$
a = \left(\frac{a_1}{p}, \frac{a_2}{p}\right) \in [0, 1]^2, \, \, 
 {\rm g.c.d.} (p, a_1, a_2) = 1, \, \, \, 
$$
$$
b=  \left(\frac{b_1}{q}, \frac{b_2}{q}\right) \in [0, 1]^2, \, \, 
 {\rm g.c.d.} (q, b_1, b_2) = 1,
 $$
 we define the operation
 $$
 a\oplus b =
 \left(\frac{a_1}{p}, \frac{a_2}{p}\right)
 \oplus
 \left(\frac{b_1}{q}, \frac{b_2}{q}\right)
       =
 \left(\frac{a_1+b_1}{p+q}, \frac{a_2+b_2}{p+q}\right).
 $$
We note that if integer vectors
 $$
 (p, a_1, a_2), \, \, (q, b_1, b_2), \, \, (r, c_1, c_2)
 $$
with corresponding points
 $$
a = \left(\frac{a_1}{p}, \frac{a_2}{p}\right), \, \, b=
\left(\frac{b_1}{q}, \frac{b_2}{q}\right), \, \, c=
\left(\frac{c_1}{r}, \frac{c_2}{r}\right), 
$$
form a basis of integer lattice, then 
for the derived points $ a\oplus b, $ and  $ a\oplus b\oplus c $ the common
 denominator and both numerators are  relatively prime 
that is $g.c.d.(p+q,a_1+b_1,a_2+b_2)=
g.c.d.(p+q+r,a_1+b_1+c_1,a_2+b_2+c_2)=1$.

 Partitions ${\rm Til}_\nu$ may be constructed as
 follows.
 The initial  partition ${\rm Til}_0$ consists of two triangles
 with vertices $(0, 0), (1, 0), (0, 1)$ and $(0, 1), (1, 0), (1, 1)$.
 Then  a triangle $\Delta$ with vertices
   $a, b, c$ in partition ${\rm Til }_\nu$
is partitioned into six triangles with vertices
  $$
  a, a\oplus b, a\oplus c;
  $$
  $$
  b, b\oplus a, b\oplus c;
  $$
  $$
  c, c\oplus a, c\oplus b;
  $$
  $$
   a\oplus b, a\oplus c, 
   a\oplus b\oplus c
   ;
  $$
  $$
  b\oplus a, b\oplus c , 
   a\oplus b\oplus c
    ;
  $$
  $$
  c\oplus a, c\oplus b , 
   a\oplus b\oplus c
  .
  $$

\includegraphics{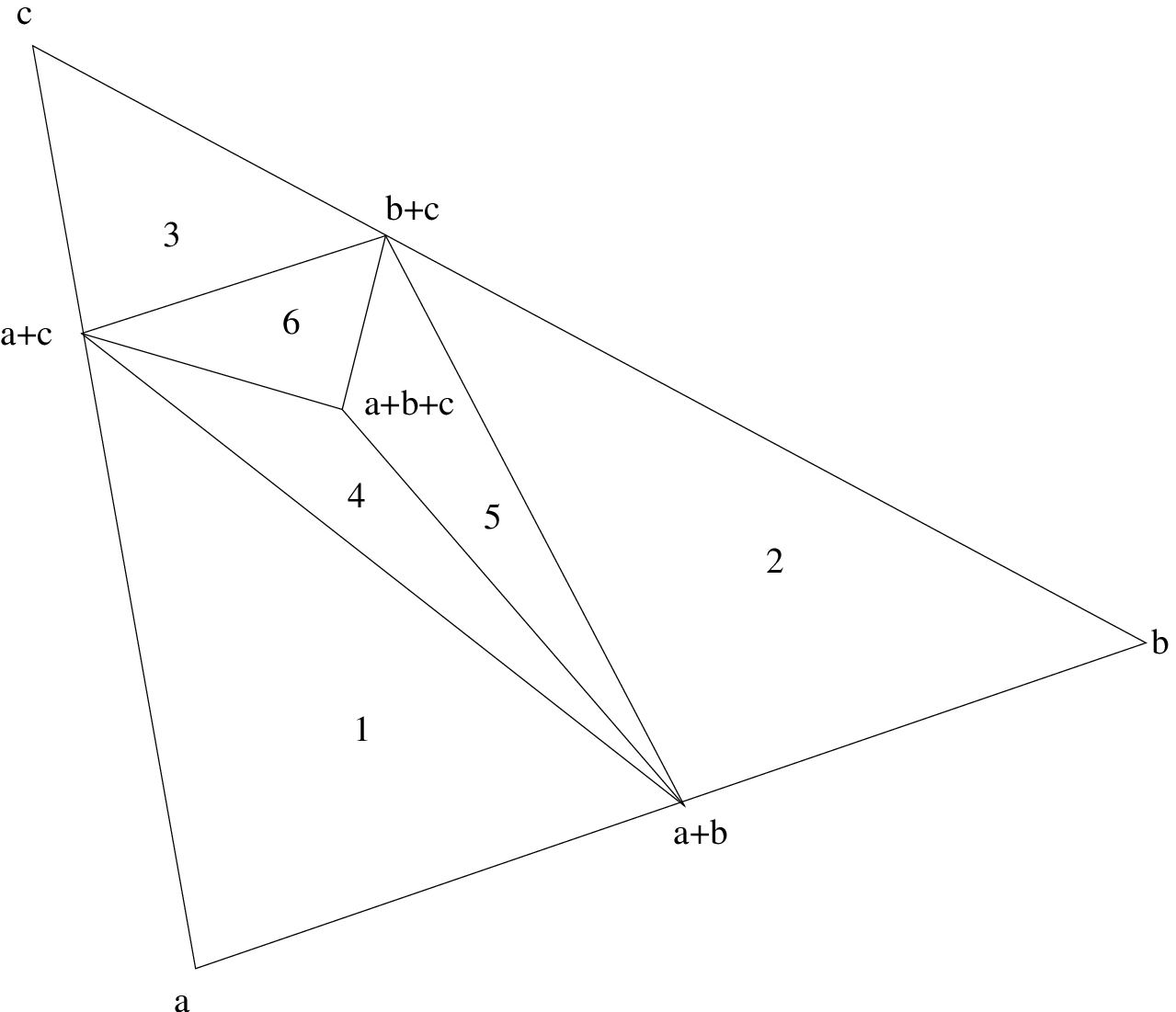}


The corresponding partition of the simplex $\Delta$ in this case  is
shown in Fig~1.

We mention a few simple combinatorial properties of the corresponding
  partitions   ${\rm Til}_\nu $ and     graphs $T_\nu , T$ and their
  respective recurrence formula:

 1. ${\rm Til}_\nu$ is a partition of the unit square $[0, 1]^2$ into
 $f_\nu = 2\times 6^{\nu}$ triangles, with $ f_{\nu} = 6f_{\nu-1}, f_0 = 2$.

 2. The number of edges of the graph $T_\nu$ is 
 $ r_\nu = 2^{\nu}\times (3^{\nu+1}+2)$, with $r_\nu = 2r_{\nu -1} +
 6f_{\nu -1}, r_0 = 5$. 

 3. The number of vertices of graph $T_\nu$ is
 $ v_\nu = 6^\nu + 2^{\nu +1} + 1$, with  
$v_\nu = v_{\nu -1} +r_{\nu - 1} + f_{\nu - 1}, v_0  = 4$.
In particular, let $v_\nu^{[d]}$ be the  the number of vertices of the
 graph $T_\nu$  with  degree $d$, then
$v_\nu^{[2]}=2$,\\
$v_\nu^{[3]}=(2\cdot 6^\nu+8)/5$, with
 $v_\nu^{[3]}=v_{\nu-1}^{[3]}+f_{\nu-1}, v_0^{[3]}=2$,\\
$v_\nu^{[5]}=2^{\nu+2}-4$, with $v_\nu^{[5]}=v_{\nu-1}^{[5]}+
 2^{\nu+1}, v_0^{[5]}=0$,\\
$v_\nu^{[8]}=(6^{\nu+1}+14)/10 - 2^{nu+1}$  
 with 
$v_\nu^{[8]}=v_{\nu-1}^{[8]}+ (3\cdot f_{\nu-1}-2^{\nu+1})/2,
 v_0^{[8]}=0$.

 4. The degree $\deg(v)$ for any vertex $v$ of the graph $T$
 takes values from the set $\{ 2, 3, 5, 8\}$.
Moreover a vertex from the set $V_\nu$ has the same degree in the
 graph $T_\nu$ and in the graph $T$. 

5. The areas of the triangles in partition Til$_\nu$ vary between
   (asymptotically) $\frac{8}{(1+\sqrt{2})^{3(\nu+1)}}$, six triangles
   obtained by 
   always using rule 6, assuming $q(c)\ge q(a),q(b)$  (use the
   recursion $x_{\nu} = 2x_{\nu-1} + x_{\nu-2}$ with initial values 1,2 for
   $q(a)_\nu, q(b)_\nu$ and 1,3 for $q(c)_\nu$), 
and (precisely) 
$\frac{1}{2(\nu+1)^2}$, the six triangles with a corner in the original
   square.

The Dirichlet series $
L(\hbox{\got A}, \beta) $ for our algorithm  can be written as follows
$$
L(\hbox{\got A}, \beta) = \sum_{a\in \mathbb{Q}^2\cap [0, 1]^2 }
\frac{\deg(a)}{ q(a)^\beta}= \sum_{q=1}^{+\infty } \frac{ 2G_2 (q)+3 
G_3(q) +5G_5(q)+8G_8(q)}{q^\beta}, 
$$
where $G_l (q), l \in \{2, 3, 5, 8\}$ denotes the number of rational
points $a\in [0, 1]^2$ with $q(a) = q$ and $\deg(a) = l.$ Clearly 
$G_2(q)+G_3(q)+G_5(q)+G_8(q)= \#\{(a_1, a_2)\in \mathbb{Z}^2\, :\, \,
\, 0\le a_1, a_2\le q, \, \, {\rm g.c.d.} (q, a_1, a_2) = 1\} \le
(q+1)^2 -4].$

$L(\hbox{\got A}, \beta)$ will be used in Section 2.5 to obtain
the asymptotic behaviour of  $\sigma_{n, \beta } (\hbox{\got A})$ for
$n \to \infty$.

%

\subsection{Lemmata about triangles from partition ${\rm Til}_\nu$}

Let a triangle $\Delta $ occur in partition ${\rm Til}_n$. We
define $\Delta^* (\Delta )$ to be the unique  triangle from 
partition 
${\rm Til }_{n-1}$ such that $ \Delta\subset \Delta^*(\Delta)$.

{\bf Lemma 4.}\, \, \, {\it Let $\Delta \in {\rm Til}_n$ 
and $a, b, c$ be the vertices of the triangle
  $\Delta^*(\Delta) $, with $a$ not a vertex of $\Delta$. Then for all
vertices $\omega$ of $\Delta $ we have   $q(\omega ) \ge
\min\{q(b),q(c)\}$.

}

\proof
The lemma is obvious, since the expression for $q(\omega)$ from
algorithm $\hbox{\got A}$ (rules 2--6) is the sum of one to three
positive summands, 
and one of them is $q(b)$ or $q(c)$. \hfill$\Box$

For every triangle $\Delta_n$ we can consider the unique sequence of
nested triangles 
\begin{equation}
\Delta_n\subset \Delta_{n-1}\subset \cdots \subset \Delta_1\subset
\Delta_0 \label{del} 
\end{equation}
 where $\Delta_\nu$ is a triangle from the partition ${\rm Til }_\nu $
 and $\Delta_\nu =\Delta^* (\Delta_{\nu +1} )$. 
Sometimes for convenience we shall write $\Delta^{[0]} (\Delta ) =
 \Delta$ and $\Delta^{[k]} (\Delta ) =   \underbrace{\Delta^* (
 ... (\Delta^*   (}_{k\, \,  \text{times}}\Delta)))$, so $\Delta_{n-k}
 = \Delta^{[k]}(\Delta_n)$. 

Now, for triangle  $\Delta$ from partition ${\rm Til }_n$  we define
 the value $t(\Delta )$ which is of principal importance for our
 proof
($t(\Delta )$ is an analogue for the partial quotient of an ordinary
 one-dimensional continued fraction).
 
If $\Delta^* (\Delta )$ has no  common 
vertices with $\Delta$ then we put $t(\Delta ) = 1$. 
If all triangles  $\Delta^{[k]}(\Delta), k = 0, ..., t$ 
 have common vertex $a$ but this vertex $a$ is not a vertex of the
 triangle $\Delta^{[t+1]}(\Delta)$, 
we write $ t(\Delta ) = t$. 
From the construction of our algorithm we
 observe that in the case $ t(\Delta ) \ge 2$ for any $k\in \{ 1, 2,
 ..., t\}$ 
the following holds: If $\Delta ^{[k]}(\Delta )$  has vertices $a, b, c$ then
 $\Delta ^{[k-1]}(\Delta )$ has vertices $a, a\oplus b, a\oplus c$.

{\bf Lemma 5.} \, \, \, 
{\it Let $\Delta $ be a triangle from the partition ${\rm Til}_n$, 
$t(\Delta ) = t$, $\Delta = \Delta^{[0]}(\Delta )\subset \cdots
\subset \Delta^{[t]}(\Delta)$, and let $a$ be the common vertex for all
these triangles.  
Then  $ a \in V_{n-t} \setminus V_{n-t-1}.$} 

\proof
From the conditions it follows that $a$ is a vertex of $ \Delta^{[t]}(\Delta)$
but not a vertex of $\Delta^*(\Delta^{[t]}(\Delta))=\Delta^{[t+1]}(\Delta)$. 
From the construction of algorithm $\hbox{\got A}$ one can see that
this may happen only when  $ a \in V_{n-t} \setminus V_{n-t-1}$.
\hfill$\Box$

After the  definition of $t(\Delta )$ we can construct a subsequence
of the sequence (\ref{del}) in the following way. Put 
\begin{equation}
\Delta_n\subset \Delta_{n-t_r}\subset \Delta_{n-t_r-t_{r-1}}\subset
 \cdots\subset\Delta_{n-t_r-t_{r-1}-...-t_2}
 \subseteq  \Delta_1 
\subset
 \Delta_{n-t_r-t_{r-1}- ...- t_2- t_1}= \Delta_0
, \label{sec}
\end{equation}
 where $t_k$ are natural numbers, 
$$
t_1+...+t_r = n
$$
and
$$
t_k = t(\Delta_{n -t_r-...-t_{k+1}}).
$$
Hence for each $\Delta$ from the partition ${\rm Til}_n$ we have the
correspondence $\Delta \mapsto [t_1, ..., t_r], \, \, t_j \in 
\mathbb{N}, \, \, t_1+...+t_r = n$. We define the sequence $[t_1, ...,
  t_r]$ as {\it code} of triangle $\Delta $ 
(different triangles from the partition ${\rm Til }_n$ may have the
same  code). 
We define the {\it empty} code to correspond to any triangle from the
initial partition ${\rm Til}_0$.

 {\bf Lemma 6.}\, \, \, {\it
Let the triangle $\Delta=\Delta_n$ with vertices $a, b, c$ and
code $[t_1, t_2, ..., t_r]$ occur in partition ${\rm Til}_n$
and $r\ge 2$. 
Let  $\Delta^{[t_r+t_{r-1}]}(\Delta)=\Delta_{n-t_r-t_{r-1}}$ 
be the triangle with code $[t_1, t_2, ..., t_{r-2}]$ and vertices
$a', b', c'$. Then 
$$
\min \{ q(a), q(b), q(c)\} \ge 2 \min \{ q(a'), q(b'), q(c')\}.
$$}

\proof
Let $\Delta^{[t_r]}(\Delta)$ have vertices $a'', b'', c''$.
We consider three cases.

1.  In the case when
 $t_r = 1$ and triangles  $\Delta , \Delta^{[t_r]}(\Delta)$ have no common
 vertex we obtain 
(for the corresponding notation of vertices)
$$
a= a''\oplus b''\oplus c'', \, \, \, b = a''\oplus c'', \, \, \, c =
a''\oplus b'', 
$$
and
$$
q(a) = q(a'')+q(b'')+q(c''), \, \, \, q(b) = q(a'')+q(c''), \, \, \,
q(c) = q(a'')+q(b''). 
$$
Hence
 $$
\min \{ q(a), q(b), q(c)\} \ge 2 \min \{ q(a''), q(b''), q(c'')\} 
\ge 2 \min \{ q(a'), q(b'), q(c')\}. 
$$

2.  In the case when
 $t_{r-1} = 1$ and triangles  $\Delta^{[t_r]}(\Delta),
 \Delta^{[t_r+t_{r-1}]}(\Delta)$ 
 have no common vertex by the same reasons we obtain 
 $$
\min \{ q(a), q(b), q(c)\} \ge \min \{ q(a''), q(b''), q(c'')\} 
\ge  2 \min \{ q(a'), q(b'), q(c')\}. 
$$

3. Finally,  we must consider the  case where the triangles  $\Delta ,
\Delta^{[t_r]}(\Delta)$ have a common vertex, and the triangles
$\Delta^{[t_r]}(\Delta),  
\Delta^{[t_r+t_{r-1}]}(\Delta)$ also have a common vertex, where the first
common vertex must not coincide with the second one. 
Without loss of generality we may 
assume that $a =a''$ is the common vertex for $\Delta ,
\Delta^{[t_r]}(\Delta)$ and $b' =b''$ is the common vertex for
$\Delta^{[t_r]}(\Delta),  
\Delta^{[t_r+t_{r-1}]}(\Delta)$. Then
$$
b =b''\underbrace{\oplus a''\cdots\oplus a''}_{t_r\, \, {\text
    times}}\, \, , c= c''\underbrace{\oplus a''\cdots\oplus
  a''}_{t_r\, \,  \text{times}},  
\, \, \, a'' =a'\underbrace{\oplus b'\cdots\oplus b'}_{t_{r-1}\, \,
   \text{times}}\, \, , c''= c'\underbrace{\oplus b'\cdots\oplus 
b'}_{t_{r-1}\, \,  \text{times}}
$$
and
$$
q(a) =q(a'') =q(a')+t_{r-1}q( b') \ge  q(a')+q( b') \ge 2 \min \{
q(a'), q(b'), q(c')\},  
$$
$$
q(b) = q(b'')+ t_r q(a'') \ge 2 \min \{ q(a'), q(b'), q(c')\}, 
$$
$$
q(c) = q(c'') + t_r q(a'') \ge 2 \min \{ q(a'), q(b'), q(c')\}.
$$
The Lemma is proved. \hfill$\Box$

 {\bf Lemma 7.}\, \, \, {\it
Let $\Delta $ with vertices $a, b, c$ have code $[t_1, ..., t_r]$. Then
$$
\min \{ q(a), q(b), q(c)\} \ge 2^{\lfloor r/2\rfloor}.
$$}

\proof
Lemma 7 follows by induction from Lemma 6. \hfill$\Box$

We need one more lemma about triangles.

{\bf Lemma 8.}\, \, \, {\it Let $\Delta$ be a triangle
from the partition ${\rm Til}_\nu $ and $a, b, c $ be vertices of $\Delta$. Let
$$
f = \min \{ q(a), q(b), q(c) \}, \, \, \, F = \max \{ q(a), q(b), q(c) \} .
$$
Then $ F\le (\nu + 1) f$. }

\proof
Induction in $\nu$. The base for $\nu = 0$ is obvious. Let $a', b', c'
$ be vertices of $\Delta^*(\Delta )$, w.l.o.g.~$q(a') 
\le q(b') \le q(c')$  and by induction assumption 
$q(c')\le \nu q(a')$. At the $\nu$-th step of algorithm
$\hbox{\got A}$, we obtain the triangle $\Delta$ from the triangle 
$\Delta^*(\Delta )$. There are six possibilities:

1. $\Delta$ has vertices $ a', a'\oplus b' , a'\oplus c'$. Then $ f =
   q(a'), F= q(a')+ q(c') \le (\nu +1 ) q(a ')$. 

2. $\Delta$ has vertices $ b', b'\oplus a' , b'\oplus c'$. Then $ f =
   q(b'), F= q(b')+ q(c') \le q( b') +  \nu  q(a ') \le (\nu +1 ) q(b
   ')$. 

3. $\Delta$ has vertices $ c', c'\oplus a' , c'\oplus b'$. Then $ f =
   q(c'), F= q(b')+ q(c') \le q( b') +  \nu  q(a ') \le (\nu +1 )
   q(c')$. 

4. $\Delta$ has vertices $ a'\oplus b' \oplus c', a'\oplus b' ,
a'\oplus c'$. Then $ f = q(a')+q(b'), F= q(a') + q(b')+ q(c')   \le
(\nu +1 ) (q(a')+q(b'))$.

5. $\Delta$ has vertices $ a'\oplus b' \oplus c', a'\oplus b' ,
b'\oplus c'$. Then $ f = q(a')+q(b'), F= q(a') + q(b')+ q(c')   \le
(\nu +1 ) (q(a')+q(b'))$.

6. $\Delta$ has vertices $ a'\oplus b' \oplus c', a'\oplus c' ,
b'\oplus c'$. Then $ f = q(a')+q(c'), F= q(a') + q(b')+ q(c')   \le
(\nu +1 ) (q(a')+q(c'))$.

So in every case we have $ F\le (\nu + 1) f$, and the lemma is
proved. \hfill$\Box$

\subsection{Global weak convergence of Algorithm {\got A}}

{\bf Theorem 1.}\, \, \, {\it 
The multidimensional continued fraction algorithm   corresponding to
algorithm $\hbox{\got A}$ weakly
  converges everywhere.}

\proof
For $\xi \in [0, 1]^2$ we look for a sequence of triangles
$$
\Delta_0\supset \Delta_1\supset \cdots \supset \Delta_\nu \supset
\cdots, \, \, \, \, \bigcap_\nu \Delta_\nu \ni \xi, \, \, \,
\Delta_\nu \in  {\rm Til }_\nu.
$$
It is sufficient to prove that ${\rm diam } \Delta_\nu \to 0 $ as $\nu
\to \infty$, where ${\rm diam} \Omega = \sup_{x, y \in \Omega } |x,
y|$,  and $|x, y|$ is the distance between points $x$ and $y$.

Let $\Delta_{\nu -1}$ have vertices
$$ a = \left(\frac{a_1}{q(a)}, \frac{a_2}{q(a)}\right), \, \, 
b= \left(\frac{b_1}{q(b)}, \frac{b_2}{q(b)}\right) , \, \, 
c=  \left(\frac{c_1}{q(c)}, \frac{c_2}{q(c)}\right)
$$

 We consider two cases.

1. Triangles $\Delta_{\nu -1}$  and $\Delta_{\nu }$ have a common
vertex. Let this common vertex be $a$. Then $\Delta_\nu $ has vertices
$ a' = a,  b' = a\oplus b, c'=a\oplus c$. Then we define the quotient
$$
s =\frac{|b',a'|}{|b,a|} = 
    \left| \frac{ \frac{ a_j+b_j}{q(a)+q(b)}-\frac{a_j}{q(a)}}{
  \frac{b_j}{q(b)}- \frac{a_j}{q(a)}} \right| =
\frac{q(b)}{q(a)+q(b)},\ \ j=1,2. 
$$
By Lemma 8 we have $ \frac{1}{\nu}\le s \le \frac{\nu -1}{\nu}$, hence
 $ \frac{1}{\nu}\le \frac{|a', b'|}{|a, b|}\le \frac{\nu -1}{\nu}$.
By the same reason  $ \frac{1}{\nu}\le \frac{|a', c'|}{|a, c|}\le
\frac{\nu -1}{\nu}$. 
So obviously 
$$
|a', b'| \le \left( 1-\frac{1}{\nu }\right) |a, b|, \, \, \, 
|a', c'| \le \left( 1-\frac{1}{\nu }\right) |a, c|. 
$$
As for $|b', c'|$ we easily deduce
$$
|b', c'| \le \left( 1-\frac{1}{\nu }\right) \max \{ |a, b|, |a, c|, |b, c|\}
$$
(if  both angles in vertices $b', c'$   are less than $\pi /2$, we have
the bound $|b, c|$; in the other cases $|a, b|$ or $|a, c|$,
respectively). 
Now 
$$
{\rm diam } \Delta_\nu = \max \{ |a', b'|, |a', c'|, |b', c'|\} \le
\left( 1-\frac{1}{\nu }\right) \max \{ |a, b|, |a, c|, |b, c|\} =$$
$$
=  \left( 1-\frac{1}{\nu }\right) {\rm diam } \Delta_{\nu - 1}.
 $$

 2. Triangles $\Delta_{\nu -1}$  and $\Delta_{\nu }$ have no common vertex.
Then $\Delta_{\nu }$  lies inside triangle $ \Delta^+$ with vertices $
a\oplus b, a\oplus c, b\oplus c$ and 
$$
{\rm diam } \Delta_\nu \le {\rm diam } \Delta^+ = \max\{ |a\oplus b,
a\oplus c|, |a\oplus b, b\oplus c|, |b\oplus c, a\oplus c| \} \le 
 \left( 1-\frac{1}{\nu }\right)
{\rm diam } \Delta_{\nu-1}
$$
by the same reasons.

In both cases, we have (using diam$(\Delta_1)=1$ for all four
$\Delta_1\in$Til$_1$)
$$
{\rm diam } \Delta_\nu \le \prod_{l = 2}^{\nu}  \left( 1-\frac{1}{l
}\right) = \frac{1}{\nu} \to 0, \, \, \, \nu \to
\infty,
$$
and the theorem is proved.\hfill$\Box$

\subsection{Asymptotic behaviour  of Algorithm {\got A}}

In this section, we will prove a formula for the moments of
$\sigma_{n,\beta}(\mbox{\got A})$ for algorithm {\got A} analogous to~(2).

 {\bf Lemma 9.}\, \, \, {\it  For any $\beta >1$, we have
$$
\sum_{n =0}^\infty
 \sigma_{n, \beta } (\hbox{\got A})
\le \frac{16}{3}\zeta (2\beta ) \zeta( 3\beta - 2 ).
$$}

\proof
For triangle $\Delta$ we consider the vertex $\alpha(\Delta )$ 
such that
 the common denominator $ q(\alpha(\Delta )) $ is the smallest among all
 vertices 
of triangle $\Delta $  (it may not be  unique and in this case  
we fix one of the minimal vertices). 
Then
 $$
\sigma_{n, \beta } (\hbox{\got A}) = \sum_{\Delta \in
  {\rm Til}_n} \left({\rm mes } \Delta \right)^\beta = \sum_{m=0}^n
\sum_{ \scriptsize
\begin{array}{c}
\Delta \in {\rm Til}_n \cr \alpha(\Delta )\in V_m\setminus V_{m-1}
\end{array}
 } \left({\rm mes } \Delta \right)^\beta
.
$$
Supposing that the following series converges (absolutely) we change
the order of summations: 
$$
\sum_{n=0}^\infty \sigma_{n, \beta } (\hbox{\got A}) =
\sum_{n=0}^\infty 
\ \ \sum_{m=0}^n 
\ \ \sum_{ 
\Delta \in {\rm Til}_n, 
 \alpha(\Delta )\in V_m\setminus V_{m-1}
 } 
\left({\rm mes } \Delta \right)^\beta =
 $$
 $$
 =
 \sum_{m=0}^\infty
\ \  \sum_{\alpha\in V_m\setminus V_{m-1}}
\ \  \sum_{n=m}^\infty
\ \   \sum_{
\Delta \in  {\rm Til}_n, \alpha(\Delta )=\alpha
 } \left({\rm mes } \Delta \right)^\beta.
$$

We fix a point $a\in V_m\setminus V_{m-1}$. 
Among the triangles from partition
${\rm Til}_m$ there are triangles with vertex $a$. (The number of
these triangles is 2, 3, 5, or 8.) 
Some of these triangles $\Delta $ may admit the
 property $a= \alpha(\Delta )$ and some may not. 
If we consider a triangle
 from ${\rm  Til}_m$ with vertices $a, b, c$ and $ a  \neq \alpha(\Delta)$
 then the  triangle $\Delta'$ 
 with vertices $ a, a\oplus b, a\oplus c$ must appear in the partition
${\rm Til}_{m+1}$, and this triangle $\Delta'$ has the property  $a=
 \alpha(\Delta ')$. 
Hence the vertex $ a\in V_m$ is totally surrounded by  triangles
\begin{equation} 
\Delta^{(1)}, ..., \Delta^{(r)}, r=\deg(a)  \in \{ 2, 3, 5, 8\},
\label{dD}
\end{equation} 
where the triangle $\Delta^{(i)}$ has vertices
$a, b^{(i)}, c^{(i)}$, $a= \alpha(\Delta^{(i)})$, and each of these
triangles 
belongs to partition ${\rm Til}_m$ or ${\rm Til}_{m+1}$. Moreover, every
triangle $\Delta $ from a partition ${\rm Til}_n, n\ge m+1$ with the 
property $a= \alpha(\Delta )$ may be obtained from one of the triangles in
(\ref{dD}) in the following sense. 
If $\Delta$ does not coincide with
one of the triangles from (\ref{dD}) and
 $a, b, c $ are vertices of $\Delta$  then
$$
b = b^{(i)}\underbrace{\oplus a\cdots\oplus a}_{j\, \,  \text{times}}
\, \, , \, \, \, c = c^{(i)}\underbrace{\oplus a\cdots\oplus a}_{j\,
  \,  \text{
times}}\, \, , 
$$
where $j = n-m$ or $ j = n-m-1$. 
We see that for the rational point $b
=  b^{(i)}\underbrace{\oplus a\cdots\oplus a}_{j\, \,  \text{times}} $
the  common denominator is $q(b) =  q(b^{(i)}) + jq(a)$ and for the
rational point $c = c^{(i)}\underbrace{\oplus a\cdots\oplus a}_{j\, \,
   \text{ times}}$ the common denominator is $q(c) = q(c^{(i)}) +
jq(a)$. 
But for any $i$ we have $ \min\{ q(a), q(b^{(i)}),
q(c^{(i)})\} = q(a) $ and  hence $ q(b), q(c) \ge (j+1)q(a)$. 
Recall that by Lemma 2,  triangle $\Delta$ with  vertices $a, b,
c$  has $ {\rm mes }\Delta = \frac{1}{2 q(a) q(b)q(c)}$.
Now we see that for fixed $a\in V_m$  we can get the upper bound
$$
  \sum_{n=m}^\infty
  \sum_{
\Delta \in {\rm Til}_n,\  
\alpha(\Delta )=a
 } \left({\rm mes } \Delta \right)^\beta
\le \frac{8}{2 (q(a))^{3\beta}} \, \, \sum_{j=1}^\infty
\frac{1}{j^{2\beta }}= \frac{4 \zeta (2\beta )}{(q(a))^{3\beta}} . 
$$
We continue our estimate:
$$
\sum_{n=0}^\infty \sigma_{n, \beta } (\hbox{\got A})
\le 4\zeta (2\beta)\times 
\left( \sum_{m=0}^\infty 
 \sum_{a\in V_m\setminus V_{m-1}}\frac{1}{(q(a))^{3\beta }}
\right).
$$
To complete the proof of  Lemma 9 we must  use the estimate
\begin{equation}
\sum_{m=0}^\infty
 \sum_{a\in V_m\setminus V_{m-1}}\frac{1}{(q(a))^{3\beta }}=
 \sum_{q=1}^\infty
 \sum_{\scriptsize
 \begin{array}{c}
 0\le a, b\le q\cr
 {\rm g.c.d.} (q, a, b) = 1
 \end{array}
 }
 \frac{1}{q^{3\beta }}
\le 
\frac{4}{3}\sum_{q=1}^\infty \frac{1}{q^{3\beta -2}}= 
\frac{4}{3}\zeta (3\beta - 2 ), \label{zeta} 
\end{equation}

where the first equality follows from the completeness property of our
algorithm, also $\ds \frac{4}{3}=\max_q\frac{(q+1)^2-4}{q^2}$. Observe
that all the series under consideration converge  (absolutely),
\hfill$\Box$

In the sequel we shall use not only Lemma 9 but also inequality
(\ref{zeta}) from the proof of Lemma 9.

Now we take parameters
\begin{equation}
\gamma = \frac{4(6\beta^2+\beta-1)}{9\log 2 \cdot (\beta - 1)\beta},
 \ \ \  
w= (\log n)^{1-\frac{1}{3\beta}} n^{\frac{2\beta +1}{3\beta}} 
 \label{par}
\end{equation}
and divide the sum from the definition $ \sigma_{n, \beta } 
(\hbox{\got  A})$ into three sums 
$$
 \sigma_{n, \beta } (\hbox{\got A})=
 \sum_{\Delta \in  {\rm Til}_n} \left({\rm mes }
 \Delta \right)^\beta = 
\Sigma_{(1)}+\Sigma_{(2)}+\Sigma_{(3)}, 
$$
where $ \Sigma_{(1)}$ is the sum over all $\Delta $ from $ {\rm
  Til}_n$ for which in the code $ [t_1, ..., t_r]$ we have 
\begin{equation}
r\ge \gamma \log n, \label{1}
\end{equation}
$ \Sigma_{(2)}$ is the sum over all $\Delta $ from $ {\rm Til}_n$ for
which in the code $ [t_1, ..., t_r]$ we have 
\begin{equation}
r< \gamma \log n, \, \, \, 1\le t_r \le n - w, \label{2}
\end{equation}
and $ \Sigma_{(3)}$ is the sum over all $\Delta $ from $ {\rm Til}_n$ for
which in the code $ [t_1, ..., t_r]$ we have 
\begin{equation}
r< \gamma \log n, \, \, \, t_r > n - w. 
 \label{3}
\end{equation}

($\Sigma_{(3)}$ will be the dominating term.)

{\bf Lemma 10.} \, \, \, 
$$
\Sigma_{(1)} \le {n^{- (3\log 2\, \, \gamma (\beta - 1))/{4}}}.
$$

\proof
Obviously,
$$
\Sigma_{(1)}\le \max_{
\Delta \in {\rm Til}_n,\ 
r\ge \gamma \log n
} 
\left({\rm mes }\Delta\right)^{\beta -1} \times 
\sum_{\Delta \in {\rm Til}_n} {\rm mes } \Delta. 
$$
Let the maximum occur on some triangle $\Delta$ with vertices $a, b, c$.
We apply (\ref{one}) and  the inequality
$$
\max_{
\Delta \in {\rm Til}_n,\ 
r\ge \gamma \log n
} \left({\rm mes }\Delta\right)^{\beta -1} = \frac{1}{(2
  q(a)q(b)q(c))^{\beta - 1}} \le \frac{1}{(2^{1+ 3\lfloor r/ 2
    \rfloor} )^{\beta - 1}}
$$
$$ 
\le \frac{1}{2^{ 3\gamma (\beta - 1) \log n/ 4}} = n^{-\frac{ 3\log
  2\, \, \gamma (\beta - 1)}{4}},
$$
which follows from Lemma 2 and Lemma 7. Lemma 10 is proved. \hfill$\Box$

{\bf Lemma 11.}
$$
\Sigma_{(2)} \le \frac{2560}{9} (\zeta (3\beta - 2))^2 \zeta (2\beta ) \left(
\frac{\gamma \log n}{w}\right)^{3\beta - 1}. 
$$

\proof
Under the conditions (\ref{2}) we see that $t_1+...+t_{r-1} > w$ and hence
  there exists $j\le  r-1 $ such that 
$ t_j \ge \tau = \left\lceil \frac{w}{\gamma  \log n-1}\right\rceil$.

For a triangle 
$\Delta$ with code $[t_1, ..., t_r]$ we consider the sequence of
  triangles (\ref{sec}) and especially the triangle $\Delta_{k}$ from
  the partition 
${\rm Til}_k$ and the next triangle $\Delta_{k+t_j}$ from the
  partition ${\rm Til}_{k+t_j}$ with $ k = t_1+...+t_{j-1}$. 
Let $a, b, c$ be the 
vertices of $\Delta_k$, then $a, b' = b\underbrace{\oplus
  a\cdots\oplus a}_{t_j\, \,  \text{times}} \, \, , \, \, \, c' = c
  \underbrace{\oplus a\cdots\oplus a}_{t_j\, \,  \text{times}}$ 
are the vertices of triangle  $\Delta_{k+t_j}$. 
For the corresponding common denominators  we have
$$
 q(b') = t_j
q(a)+q(b) \ge t_j q(a) , \, \, q(c') = t_j q(a) + q(c) \ge t_j q(a).
$$
As  the element $t_j$ is {\it not the last} element in the code $[t_1,
  ..., t_r]$ in the complete sequence of triangles~(\ref{del}), there
exists the triangle $\Delta_{k+t_j + 1}$. 
By Lemma 4, for every
vertex $\omega $ of the triangle $\Delta_{k+t_j + 1}$ we have 
\begin{equation}
q(\omega )\ge t_j q(a).
\label{from4}
\end{equation}
Now we look for the partition ${\rm Til}_n$ restricted to the triangle
$\Delta_{k+t_j + 1}$. It is isomorphic to the  
partition $ {\rm Til}_{n - k-t_j -1}$. 
Moreover  for any triangle $\Delta \subset \Delta_{k+t_j + 1}$ with
vertices $s, u, v$ from the partition 
${\rm Til}_n$ and the corresponding triangle $\Delta '$ with vertices
$s', u', v'$ from the isomorphic partition $ {\rm Til}_{n - k-t_j -1}$,
by  (\ref{from4}) we deduce that $ q(s) \ge t_jq(a)\cdot q(s'), q(u)
\ge t_jq(a)\cdot q(u'), q(v) \ge t_jq(a)\cdot q(v')$, and hence 
$$
{\rm mes}\Delta = \frac{1}{2q(s)q(u)q(v) }\le
\frac{1}{2(t_jq(a))^3q(s')q(u')q(v') } = \frac{{\rm mes }\Delta
  '}{(t_jq(a))^3}  . 
$$

On the other hand,  Lemma 5 shows that  the vertex $a$ of the triangle
$\Delta_k$ satisfies $ a \in V_k\setminus V_{k-1}$. We take into
account  that vertex $a$ may be a common vertex for no more than 
eight triangles from the partition $ {\rm Til}_k$. 
Also we must take into account that in partition ${\rm Til}_n$ there exist just
five triangles $\Delta$ satisfying the conditions of  Lemma 4 with the
given  $\Delta^*(\Delta )$. Hence
$$
\Sigma_{(2)} 
\le 
\sum_{\tau \le t\le n} 
\sum_{\scriptsize
\begin{array}{c}
k\ge 0, h\ge 1: \cr k+h = n-t
\end{array}
}  
\left(
\sum_{a\in V_k\setminus V_{k-1}}
\frac{8}{(q(a))^{3\beta}}\right)\times \frac{1}{t^{3\beta }} \times
\left( \sum_{\Delta \in {\rm Til}_{h-1} } 5 \left({\rm mes
}\Delta\right)^{\beta} \right)\le 
$$
$$\le
\frac{40}{\tau^{3\beta - 1}} \, \, \times \, \, \left(\sum_{k=0}^\infty
\sum_{a\in V_k\setminus V_{k-1}}
\frac{1}{(q(a))^{3\beta}}\right) \, \,  
\times \, \,  \left(\sum_{h=0}^\infty \sum_{\Delta \in {\rm Til}_{h} }
\left({\rm mes }\Delta\right)^{ \beta} 
\right). 
$$

But
$$
\sum_{h=0}^\infty  \sum_{\Delta \in {\rm  Til}_{h} }  
\left({\rm mes }\Delta\right)^{\beta}  =
\sum_{h=0}^\infty \sigma_{h, \beta } (\hbox{\got A}) 
\le \frac{16}{3}\zeta(2\beta)\zeta(3\beta-2) 
$$
by Lemma 9 and
$$
\sum_{k=0}^\infty  \sum_{a\in V_k\setminus V_{k-1}}
\frac{1}{(q(a))^{3\beta}}
\le \frac{4}{3} \zeta(3\beta-2),  
$$
applying the upper bound (\ref{zeta}). Now the inequality of Lemma 11
follows.\hfill$\Box$

{\bf Lemma 12.}
$$
\Sigma_{(3)} = \frac{L(\hbox{\got A}, 3\beta)}{(2n^2)^{\beta} }+
O\left(\Sigma_{(1)}\right)+ O\left(\frac{1}{n^{2\beta} w^{3(\beta - 
1)}}+\frac{w}{n^{2\beta + 1}}\right)
$$

\proof
Obviously
$$
\Sigma_{(3)} = \Sigma_{(3)}' + O\left(\Sigma_{(1)}\right),
$$
where
$$
\Sigma_{(3)}' = \sum_{\scriptsize
\begin{array}{c}
\Delta \in {\rm Til}_n:\cr \text{code}\, \,
       [t_1, ..., t_r] 
\, \, \text{such}\, \, \text{that }\, \,  t_r > n - w
\end{array}}
\left({\rm mes }\Delta\right)^\beta.
$$

 Let $\Delta $ be a triangle from ${\rm Til}_n$ with code
$[t_1, ..., t_{r-1}, t_r]$ and $ t_r > n -w$. Then the
 triangle $\Delta ' = 
\Delta^{[t_r]}$ belongs to the partition ${\rm Til}_{n- t_r}$ and its code is
$[t_1, ..., t_{r-1}], t_1+...+t_{r-1} = n - t_r < w$.

Define $a$ to be the
common vertex for $\Delta , \Delta^{[t_r]}$. Then by Lemma 5 we have $a \in
V_{n-t_r}\setminus V_{n-t_r-1}$.

On the other hand, for any triangle $\Delta '$ 
with code $[t_1, ...,  t_{r-1}]$
from partition ${\rm Til}_{m}$, $m=n-t_r<w$ with fixed vertex 
$ a \in  V_{m}\setminus V_{m-1}$, there exists only one triangle in
Til$_n$ with code $[t_1, ..., t_{r-1}, t_r]$ and   vertex $a$. 
Hence
$$
\Sigma_{(3)}' = \sum_{m=0}^{w-1} \sum_{ a \in V_{m}\setminus V_{m-1}}
\sum_{
\scriptsize
\begin{array}{c}
\Delta \in {\rm Til}_m:  \cr a \, \, 
\text{is}\, \,  \text{
a }\, \,  \text{vertex}\, \,  \text{of} \, \,  \Delta
\end{array}}
\left({\rm mes }\Delta\right)^\beta, 
$$
where $\Delta\subset \Delta '$ is the unique triangle with code
$[t_1, ..., t_{r-1}, t_r]$ 
and common vertex $a$.
Let $\Delta '$ from ${\rm Til}_m$ have vertices $a, b, c.$
Then $\Delta $ has vertices $ a,$ $b\underbrace{\oplus a\cdots\oplus
a}_{t_r\, \,  \text{times}},$ 
$c\underbrace{\oplus a\cdots\oplus a}_{t_r\, \, \text{times}}, $  
and by Lemma 8
$$
q( b\underbrace{\oplus a\cdots\oplus a}_{t_r\, \,  \text{times}})
= t_rq(a)+ q(b)\le (n+1)q(a), 
$$
$$
q( c\underbrace{\oplus a\cdots\oplus a}_{t_r\, \,  \text{times}})
= t_rq(a)+ q(c)\le (n+1) q(a).
$$
Recall that $ t_r = n-m> n - w$. 
Applying Lemma 2, we have
$$
\sum_{m=0}^{w-1} 
\sum_{ a \in V_{m}\setminus V_{m-1}} 
\frac{\deg(a)}{(2(n+1)^2q(a)^3  )^\beta}
 \le
\Sigma_{(3)}' 
\le \sum_{m=0}^{w-1} \sum_{ a \in V_{m}\setminus V_{m-1}} 
\frac{\deg(a)}{(2(n-m)^2q(a)^3)^\beta},
$$
and thus
$$
\Sigma_{(3)}' = \sum_{m=0}^{w-1} \sum_{ a \in V_{m}\setminus V_{m-1}}
\frac{\deg(a)}{(2n^2q(a)^3)^\beta} \left( 
1+O\left(\frac{w}{n}\right)\right).
$$

But
\begin{equation}
\sum_{q=1}^{w-1} \frac{\sum_l lG_l(q)}{q^{3\beta}} \le
 \sum_{ a \in V_{w-1}} \frac{\deg(a)}{(q(a))^{3\beta}}
 \le
\sum_{q=1}^\infty \frac{\sum_l lG_l(q)}{q^{3\beta}}
\label{SumDeg1}
\end{equation}
and
$$
\sum_{q=w}^\infty \frac{\sum_l lG_l(q)}{q^{3\beta}} \ll \sum_{q=w}^\infty
\frac{1}{q^{3\beta-2}} \ll_\beta w^{-3(\beta - 1)}.
$$
Hence
\begin{equation}
 \sum_{m=0}^{w-1} \sum_{ a \in V_{m}\setminus V_{m-1}}
 \frac{\deg(a)}{(q(a))^{3\beta}} 
= \sum_{ a \in V_{w-1}}  \frac{\deg(a)}{(q(a))^{3\beta}} 
= L(\hbox{\got A}, 3\beta) +O\left(w^{-3(\beta - 1)} 
\right).
\label{SumDeg2}
\end{equation}

It follows that
$$
\Sigma_{(3)}' = \left(
 \frac{L(\hbox{\got A}, 3\beta)}
{(2n^2)^\beta } + O\left( \frac{1}{n^{2\beta }w^{3(\beta - 1)}}
\right)\right) \left( 1+O\left( \frac{w}{n}\right)\right), 
$$
and the lemma is proved.\hfill $\Box$

{\bf Theorem 2.}\, \, \, \, {\it For $\beta > 1$ the following
  asymptotic formula is valid 
$$
\sigma_{n, \beta } (\hbox{\got A}) = \frac{L(\hbox{\got A},
  3\beta)}{(2 n^2)^{\beta} }\left(1+ O\left( \frac{(\log 
n)^{1-\frac{1}{3\beta}}}{n^{\frac{\beta -1}{3\beta}}}\right)\right) .
$$
}

\proof
We need to put together the results of Lemmata 10, 11, 12 and take into
account the choice of parameters in (\ref{par}):
$$
\sigma_{n, \beta } (\hbox{\got A}) = \frac{L(\hbox{\got A},
  3\beta)}{(2 n^2)^{\beta} }+ O\left( \frac{1}{n^{\frac{ 3\log 2\, \,
      \gamma (\beta - 
1)}{4}}} + \frac{1}{n^{2\beta} w^{3(\beta - 1)}}+\frac{w}{n^{2\beta +
    1}} + \left( \frac{\gamma \log n}{w}\right)^{3\beta - 1} 
 \right)=
 $$
 $$
 = \frac{L(\hbox{\got A}, 3\beta)}{(2
n^2)^{\beta} }\left(1+ O\left( \frac{(\log
n)^{1-\frac{1}{3\beta}}}{n^{\frac{\beta -1}{3\beta}}}\right)\right).
 $$
This shows the asymptotic formula.\hfill$\Box$

\section{Algorithm $\hbox{\gotbig B}$}

\subsection{The description of Algorithm $\hbox{\got B}$}

We fix the initial partition of the unit square $\{z=(x, y_1, y_2)\,
 :\, \, x=1, y_{1, 2} \in [0, 1]\}$ into two triangles 
 $\Delta^{0,1}$ with vertices $(1, 0, 0), (1, 1, 0), (1, 0, 1)  $
and $\Delta^{0,2}$ with vertices $ (1, 1, 1), (1, 0, 1), (1, 1, 0)$.  
The vertices of both triangles form bases ${\cal E}^{0,1}$ 
and ${\cal E}^{0,2}=\{(1,1,1),(1,0,1),(1,1,0)\}$  (the {\it order} is
of importance here) of the integer lattice $\mathbb{Z}^3$. 
Now we suppose that a basis
$$
{\cal E}^{\nu,j} = \{ g^{\nu , j}_{1}, g^{\nu , j}_2, g^{\nu , j}_3\}
$$
(and corresponding triangle $ {\Delta}^\nu_j $ of the partition of the
unit square $[0, 1]^2$ into triangles) 
 occurs in our algorithm and  we define
the rule for constructing the bases for the next step of algorithm
(the rule for dividing cone ${\cal C}({\cal E}^{\nu,j})$ of the basis
${\cal E}^{\nu,j}$).
 In our algorithm $\hbox{\got B}$ the rule will be the same for
each step of the algorithm and for each basis. Namely, for the basis
${\cal E}^{\nu,j}$ which occurs at $\nu$-th step 
 we take $2$ bases ${\cal E}^{\nu +1,2(j-1)+i}, \, i = 1,2$
by the following formulas
$$
{\cal E}^{\nu+1,2(j-1)+1} = \{g^{\nu , j}_{2}+g^{\nu, j}_3, 
g^{\nu,j}_1, g^{\nu , j}_2\},  \ \ \ \text{operation\ ``1''}
$$
$$ 
{\cal E}^{\nu+1,2(j-1)+2} = \{g^{\nu , j}_{2}+g^{\nu, j}_3, 
g^{\nu , j}_1, g^{\nu , j}_3\}  \ \ \ \text{operation\ ``0''}.
$$
The construction of the set of bases ${\cal E}^{\nu +1,2(j-1)+i}, \,
1\le i\le 2$ depends on the {\it order} of elements of the  
basis~${\cal E}^{\nu ,j}$.

Obviously, this rule satisfies the conditions \ref{cond} $(i), (ii)$
--- each new set of vectors ${\cal E}^{\nu +1,2(j-1)+i}$
is a basis of the integer lattice, and the cones
${\cal C}({\cal E}^{\nu+1,2(j-1)+i}),$ $  \, 1\le i\le 2$
form a regular partition of the cone ${\cal C}({\cal E}^{\nu,j})$. 
Clearly  algorithm $\hbox{\got B}$ is finite and hence complete.
The construction described was introduced and studied in \cite{MON}
(for the general $d$-dimensional situation). For example, in \cite{MON} it
was shown that the corresponding multidimensional continued fraction
algorithm weakly converges.

\subsection{Algorithm $\hbox{\got B}$ in terms of
  constructing rational   points in the square $[0, 1]^2$}

Note that if integer vectors
$(p, a_1, a_2), \, \, (q, b_1, b_2)$
can be  extended to a basis
of the integer lattice $\mathbb{Z}^3$, then the corresponding  rational points
$a = \left(\frac{a_1}{p}, \frac{a_2}{p}\right)$ and $b=
\left(\frac{b_1}{q}, \frac{b_2}{q}\right)  $ define a rational point 
$a\oplus b $, whose common denominator and both numerators are relatively
prime. 

Partitions ${\rm Til}_\nu$ may be constructed as
follows.
The initial  partition ${\rm Til}_0$ consists of two triangles
with vertices $(0, 0), (1, 0), (0, 1)$ and $(1, 1), (1, 0), (0, 1)$.
Then  a triangle $\Delta$ with vertices
$a, b, c$ in partition ${\rm Til }_\nu$
must be partitioned into two triangles with vertices
$$
c\oplus b, a, b, \, \, \,  \text{and }\, \, \,  c\oplus b , a, c.$$ 
We must note that the {\it order} of the enumeration of vertices of
triangle $\Delta$ is important for the constructing of our partition. 

We call the first rule an operation ``1'', the second one an operation
``0''. To every triangle $\Delta\in \text{Til}_n$, we then attach a
code $\co(\Delta) = \co_1\dots \co_n$, where $\co_k\in\{0,1\}$ states, which
rule was used for the $k$-th partition. Also, let
$|\co|(\Delta)=\sum_{k=1}^n \co_k$ be the number of operations ``1'' to
obtain $\Delta$. 


\includegraphics{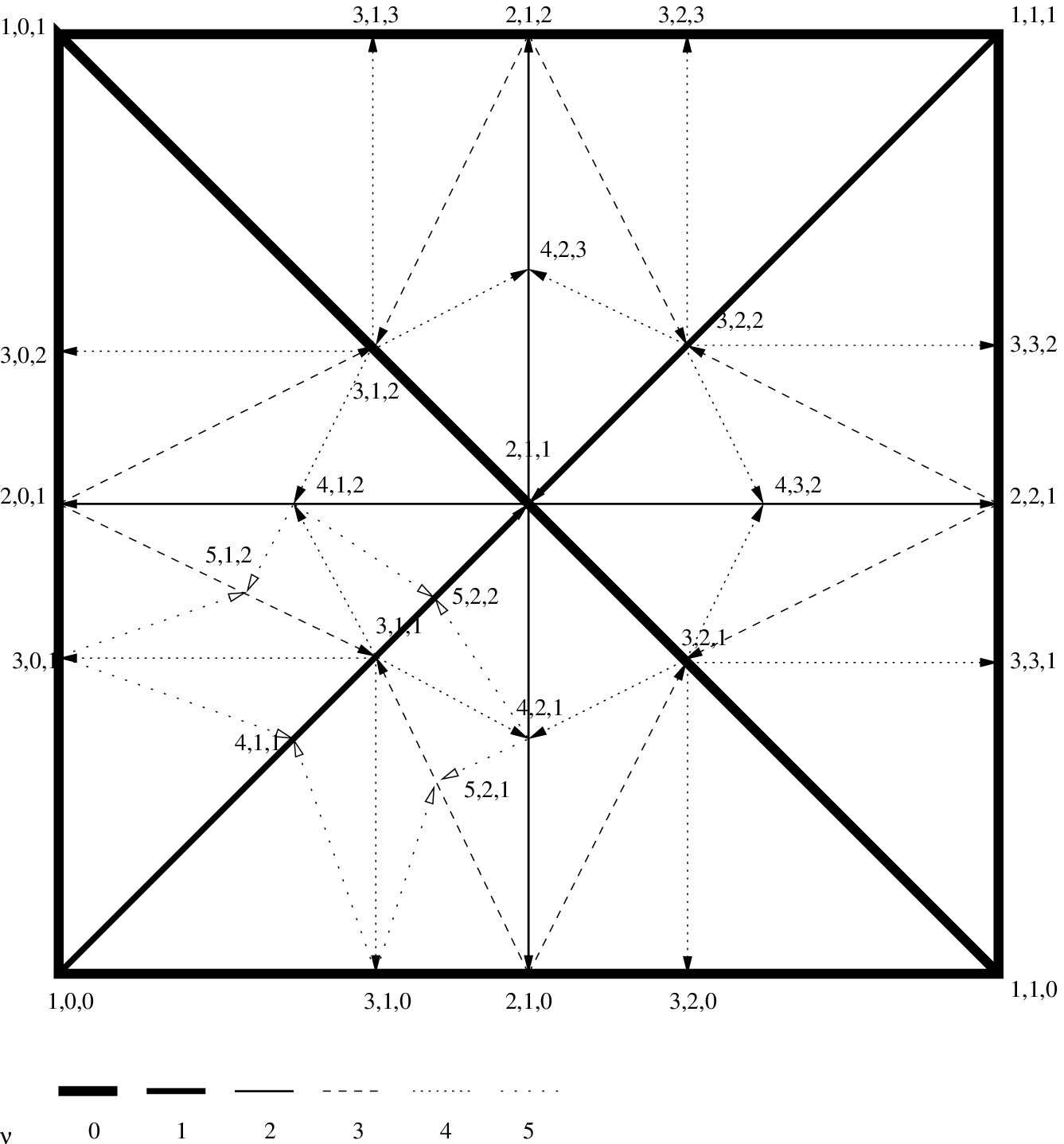}

\vspace*{5 mm}

Figure 2 shows the first 6 partitions (Til$_5$ only in the lower left
quarter), together with the points $(q,a_1,a_2)$. Walking in the
direction of the arrowheads, one obtains   a binary tree, infinite
in the case of $T$.

Lemmata 2 and  8 (consider 
$q(b)+q(c)\le \nu q(c)+q(c)=(\nu+1)q(c)\le (\nu+1)q(b)$ and 
$q(b)+q(c)\le \nu (\nu+1)q(c)$ in case of operations ``0'', and ``1'',
respectively) and inequality  (\ref{zeta}) remain valid for algorithm {\got B}.

{\bf Lemma 13.}\, \, \, \\
{\it 
$(i)$ We have $q(b)+q(c)\ge q(a)\ge q(b)\ge q(c)$ for all triangles.

\noindent$(ii)$ Let the triangle $\Delta'$ be obtained by rule ``$1$''
from triangle $\Delta$. Then 
$\text{mes}\Delta'\le \frac{1}{2}\text{mes}\Delta$.

\noindent$(iii)$ Let the triangle $\Delta'$ with vertices $a',b',c'$  
be obtained by applying rule ``$0$'' $k$ times to triangle $\Delta$ with
vertices $a,b,c$. 
Then the vertices of $\Delta'$ are
$(a',b',c') =
(a\underbrace{\oplus c\cdots\oplus c}_{k/2\ \text{times}}, 
 b\underbrace{\oplus c\cdots\oplus c}_{k/2\ \text{times}},
c)$, if $k$ is even, and
$(a',b',c') =
(b\underbrace{\oplus c\cdots\oplus c}_{(k+1)/2\ \text{times}}, 
 a\underbrace{\oplus c\cdots\oplus c}_{(k-1)/2\ \text{times}},
c)$, if $k$ is odd.
Also, $q(a'),q(b') \ge \frac{k+1}{2} q(c)$. 
}

\proof
$(i)$ By induction, obvious for $\nu=0$.
Let $(a',b',c')$ be the new vertices. Then\\
$q(b')+q(c') = q(a)+q(b)\ge q(a')=q(b)+q(c)\ge q(b') = q(a) \ge
q(c')=q(b)$ and\\
$q(b')+q(c') = q(a)+q(c)\ge q(a')=q(b)+q(c)\ge q(b') = q(a) \ge
q(c')=q(c)$, respectively. 

$(ii)$ From $(i)$, $q(b)+q(c)\geq 2\cdot q(c)$, and using Lemma 2,
$$\text{mes} \Delta'  = (2(q(b)+q(c))\cdot q(a)\cdot q(b))^{-1}
\le   (2(q(a)\cdot q(b)\cdot 2q(c))^{-1}=\frac{1}{2}\text{mes} \Delta .$$

$(iii)$ Observe the effect of rule ``0'' on triangle $\Delta$:
$(a,b,c)
\to(b\oplus c,a,c)
\to(a\oplus c,b\oplus c,c)
\to(b\oplus c\oplus c,a\oplus c,c)
\to(a\oplus c\oplus c,b\oplus c,c)
\to\dots$ Iterating the rule $k$ times leads to the stated formula,
and then 
$
q(a') = q(a)+\frac{k}{2} q(c), 
q(b') = q(b)+\frac{k}{2}q(c), 
q(c') = q(c)$, if $k$ is even, and
$
q(a') = q(b)+\frac{k+1}{2} q(c), 
q(b') = q(a)+\frac{k-1}{2}q(c), 
q(c') = q(c)$, if $k$  is odd, in any case $q(a'), q(b') \ge 
(\frac{k-1}{2}+1) q(c).$\hfill$\Box$

It is easy to verify the following properties of partitions
 ${\rm Til}_\nu $ and     graphs $T_\nu , T$ by induction.

1. ${\rm Til}_\nu$ is a partition of the unit square $[0, 1]^2$ into
$f_\nu =  2^{\nu+1}$ triangles.

2. The number of edges of the graph $T_\nu$ is equal to
$$
 r_\nu =
\begin{cases}
3\times 2^{2k}+2^{k+1}, \, \, \,  \text{if}\, \, \, n = 2k, \cr
6\times 2^{2k}+2^{k+1}, \, \, \,  \text{if}\, \, \, n = 2k+1. 
\end{cases}
$$

3. The number of vertices of graph $T_\nu$ is equal to
$$
 v_\nu =
\begin{cases}
(2^k+1)^2, \, \, \,  \text{if}\, \, \, n = 2k, \cr (2^{k}+ 1)^2
  +2^{2k}, \, \, \,  \text{if}\, \, \, n = 2k+1. 
\end{cases}
$$

4. The degree $\deg(v)$ for any vertex $v$ of the graph $T$
takes values from the set $\{3, 5, 8\}$.
In each graph $T_\nu$ also occur vertices of degree  2 (for $\nu=0$)
or 4 (for $\nu\ge 1$); 
these vertices lie in $V_\nu \setminus V_{\nu - 1}$.
The number of vertices from $T_\nu$ with the given degree can be easily
calculated.

 The Dirichlet series $
L(\hbox{\got B}, \beta) $ for our algorithm  can be written as follows
$$
L(\hbox{\got B}, \beta) = \sum_{a\in  \mathbb{Q}^2\cap [0, 1]^2 }
\frac{\deg(a)}{ q(a)^\beta}= \sum_{q=1}^{+\infty } \frac{ 3 G_3(q) 
+5G_5(q)+8G_8(q)}{q^\beta}, 
$$
where $G_l (q), l \in \{3, 5, 8\}$ is the number of rational points
$a\in [0, 1]^2$ with common denominator $q(a) = q$ and $\deg(a) = l.$ 
Obviously $G_3(q)+G_5(q)+G_8(q)= \#\{(a_1, a_2)\in \mathbb{Z}^2\, :\,
\, \, 0\le a_1, a_2\le q, \, \,$ 
${\rm g.c.d.} (q, a_1, a_2) = 1\} \le (q+1)^2.$

{\bf Lemma 14.}\, \, \, 
{\it
For all $\beta>1$, we have
$$
\sum_{n =0}^\infty
 \sigma_{n, \beta } (\hbox{\got B})
\le 
\frac{32}{3}2^{\beta}\zeta (2\beta ) \zeta( 3\beta - 2 ).
$$
}

\proof
We follow closely the proof of Lemma 9. In algorithm {\got B},
always $\alpha(\Delta) = c$ by Lemma~13(i). Hence
$$
\sum_{n=0}^\infty \sigma_{n, \beta } (\hbox{\got B}) =
\sum_{m=0}^\infty
\ \  \sum_{c\in V_m\setminus V_{m-1}}
\ \  \sum_{n=m}^\infty
\ \   \sum_{
\Delta \in  {\rm Til}_n, \alpha(\Delta )=c
 } \left({\rm mes } \Delta \right)^\beta.
$$

We fix a point $c\in V_m\backslash V_{m-1}$. 
Then we have triangles 
\begin{equation} 
\Delta^{(1)}, ..., \Delta^{(\deg(c))}\in\text{Til}_m\cup\text{Til}_{m+1},
\label{dD3}
\end{equation} 
all including vertex $c$.
Again, every
triangle $\Delta \in {\rm Til}_n$ with $\alpha(\Delta)=c$ 
is included in some $\Delta^{(i)}$
with vertices $a',b',c'=c$, and has been obtained by operations
``0''. 
With Lemma 2 and Lemma 13(iii),
$$
  \sum_{n=m}^\infty
  \sum_{\Delta \in {\rm Til}_n,\  \alpha(\Delta ) = c} 
\left({\rm mes } \Delta \right)^\beta
\le 
\frac{8}{(2 q(c)^{3})^\beta} \times 
\sum_{j=1}^\infty
\frac{2^{2\beta}}{(j+1)^{2\beta }}
\le \frac{2^{3+\beta}}{q(c)^{3\beta}}\zeta (2\beta ),
$$
and thus
$$
\sum_{n=0}^\infty \sigma_{n, \beta } (\hbox{\got B})
\le 2^{3+\beta} \zeta (2\beta)\times 
\sum_{m=0}^\infty 
\sum_{c\in V_m\setminus V_{m-1}}\frac{1}{q(c)^{3\beta}}.
$$
Using (\ref{zeta}), we obtain the result as in the proof of 
Lemma 9.\hfill$\Box$ 

We choose parameters
$$\gamma = \frac{4(6\beta^2+\beta-1)}{3\log 2 \cdot (\beta - 1)\beta},
\ \ \text{and}\ \ 
w= (\log n)^{1-\frac{1}{3\beta}}
n^{\frac{2\beta +1}{3\beta}}$$ 
and again we divide $\sigma_{n, \beta } 
(\hbox{\got  B})$ into three sums, now according to $\co$ and $|\co|$,
$$
 \sigma_{n, \beta } (\hbox{\got B})=
 \sum_{\Delta \in  {\rm Til}_n} \left({\rm mes }
 \Delta \right)^\beta = 
\Sigma_{(1)}+\Sigma_{(2)}+\Sigma_{(3)}, 
$$
where $ \Sigma_{(1)}$ is the sum over all $\Delta $ from 
$ {\rm  Til}_n$ with
\begin{equation}
|\co|(\Delta)\ge \gamma \log n, \label{1b}
\end{equation}
$ \Sigma_{(2)}$ is the sum over all $\Delta $ from 
$ {\rm Til}_n$ with
\begin{equation}
|\co|(\Delta)< \gamma \log n, \, \, \, \exists k>w\colon \co_k=1, \label{2b}
\end{equation}
and $ \Sigma_{(3)}$ is the sum over all $\Delta $ from 
$ {\rm Til}_n$ with
\begin{equation}
|\co|(\Delta)< \gamma \log n, \, \, \, \co_{w+1}=\cdots=\co_n=0.
 \label{3b}
\end{equation}

{\bf Lemma 15.}\, \, \, For all $\beta > 1$,
$$
\Sigma_{(1)} \le {n^{- \log 2\, \, \gamma (\beta - 1)}}.
$$

\proof
Obviously,
$$
\Sigma_{(1)}\le \max_{
\Delta \in {\rm Til}_n,\ 
|\co|(\Delta)\ge \gamma \log n
} 
\left({\rm mes }\Delta\right)^{\beta -1} \times 
\sum_{\Delta \in {\rm Til}_n} {\rm mes } \Delta. 
$$
Let the maximum occur on some triangle $\Delta$ with vertices $a, b, c$.
We apply Lemma 13$(ii)$, (\ref{one}) and  the inequality
$$
\max_{
\Delta \in {\rm Til}_n,\ |\co|(\Delta)\ge \gamma \log n
} 
\left({\rm mes }\Delta\right)^{\beta -1} 
\le 
\frac{1}{(2^{ |\co|(\Delta)})^{\beta - 1}} 
\le 
\frac{1}{2^{ \gamma (\beta - 1) \log n}} 
= 
n^{-\log 2\, \, \gamma (\beta - 1)}.
$$
Lemma 15 is proved. \hfill$\Box$

{\bf Lemma 16.} \, \, \, 
{\it 
Let some triangle $\Delta=(a,b,c)$ be given,
$\Delta'=(a',b',c')$ resulting from $\Delta$ via the operation
$\delta_0\in\{0,1\}$, and $\Delta''=(a'',b'',c'')$ resulting from $\Delta'$ via
the three operations $\delta_1,1,0$ with $\delta_1\in\{0,1\}$.
Then $c''$ was not a vertex of $\Delta, c''\not\in\{a,b,c\}$,
but it is a vertex of $\Delta'$.
}

\proof
We verify the four cases:

$(a,b,c)
\stackrel{0}{\to}(b\oplus c,a,c)
\stackrel{0}{\to}(a\oplus c,b\oplus c,c)
\stackrel{1}{\to}(b\oplus c\oplus c, a\oplus c, b\oplus c)
\stackrel{0}{\to}(a\oplus b\oplus c\oplus c, b\oplus c\oplus c,b\oplus c)
$ with $c''=b\oplus c=a'$.

$(a,b,c)
\stackrel{0}{\to}(b\oplus c,a,c)
\stackrel{1}{\to}(a\oplus c,b\oplus c,a)
\stackrel{1}{\to}(a\oplus b\oplus c, a\oplus c, b\oplus c)
\stackrel{0}{\to}(a\oplus b\oplus c\oplus c, a\oplus b\oplus c,b\oplus c)
$ with $c''=b\oplus c=b'$.

$(a,b,c)
\stackrel{1}{\to}(b\oplus c,a,b)
\stackrel{0}{\to}(a\oplus b,b\oplus c,b)
\stackrel{1}{\to}(b\oplus b\oplus c, a\oplus b, b\oplus c)
\stackrel{0}{\to}(a\oplus b\oplus b\oplus c, b\oplus b\oplus c,b\oplus c)
$ with $c''=b\oplus c=b'$.

$(a,b,c)
\stackrel{1}{\to}(b\oplus c,a,b)
\stackrel{1}{\to}(a\oplus b,b\oplus c,a)
\stackrel{1}{\to}(a\oplus b\oplus c, a\oplus b, b\oplus c)
\stackrel{0}{\to}(a\oplus b\oplus b\oplus c, a\oplus b\oplus c,b\oplus c)
$ with $c''=b\oplus c=b'$.

Observe that in any case, $c''=b\oplus c\not\in\{a,b,c\}$.\hfill$\Box$

{\bf Lemma 17.}
$$
\Sigma_{(2)} \le \frac{2^{8}\cdot 59}{9}\cdot 432^\beta  
(\zeta (3\beta - 2))^2 \zeta (2\beta)
\left(\frac{\gamma \log n}{w}\right)^{3\beta - 1}. 
$$

\proof
Condition (\ref{2b}) implies $|\co|(\Delta) < \gamma  \log n$.
As the last operation ``1'' occurs after the first $w$ partitions,
for some  $t\ge \tau = \left\lceil \frac{w}{\gamma  \log
  n}\right\rceil-1$,  
there exists $k\le w$ such that 
$(i)$ $\co_{k+1}=\dots=\co_{k+t}=0$ ($t$ consecutive operations ``0''), 
$(ii)$ $\co_{k+t+1}=1$,
$(iii)$ $\co_k=1$, $\co_{k-1}=\delta\in\{0,1\}$, or $0\le k\le 1$.

We consider first the case $0\le k\le 1$:
This part adds at most $8\cdot(\#V_0+\#V_1)$ summands of the form
$q(c)^{-3\beta}\le 1$. Since $\#V_0=4, \#V_1=5$, this amounts to at most 
$\sum_{k=0,1}\le 72$.

Let now $k\ge 2$.
For the triangle $\Delta$ with code $[\co_1, ..., \co_n]$ we consider 
the sequence of  triangles (\ref{del}), especially the triangle
  $\Delta_{k+1}\in {\rm Til}_{k+1}$ and the triangle
$\Delta_{k+t}\in{\rm Til}_{k+t}$. 
Let $\Delta_{k+1} =(a,  b, c)$ and $\Delta_{k+t} = (a',b',c')$.

By Lemma 13(iii), the corresponding  common denominators satisfy 
$$ 
q(c') = q(c),\ \ 
q(a') \ge \left(\left\lfloor\frac{t}{2}\right\rfloor-1\right) q(c) , \ \  
q(b') \ge \left(\left\lfloor\frac{t}{2}\right\rfloor-1\right) q(c).
$$
Now, the triangle $\Delta_{k+t+1}$ is obtained from $\Delta_{k+t}$ by an
operation ``1'' ($\co_{k+t+1}=1$), 
hence for every vertex $\omega\in\{b'\oplus c',a',b'\} $ of the triangle 
$\Delta_{k+t + 1}$, we have 
\begin{equation}
q(\omega )\ge \left(\left\lfloor\frac{t}{2}\right\rfloor-1\right) q(c)
\ge\frac{t}{3}\cdot q(c).
\label{from4B}
\end{equation}
Now we consider the partition ${\rm Til}_n$ restricted to the triangle
$\Delta_{k+t + 1}$. It is isomorphic to the  
partition $ {\rm Til}_{n - k-t -1}$. 
Moreover,  for any triangle
$\Delta \subset \Delta_{k+t + 1}$ with vertices $s, u, v$ from the
partition 
${\rm Til}_n$ and the corresponding triangle $\Delta '$ with vertices
$s', u', v'$ from the isomorphic partition $ {\rm Til}_{n - k-t -1}$,
by  (\ref{from4B}), we deduce that 
$$
{\rm mes}\Delta = \frac{1}{2q(s)q(u)q(v) }
\le \frac{1}{2(\frac{t}{3}\cdot q(c))^3 q(s')q(u')q(v') } 
= \frac{27{\rm mes }\Delta '}{(t q(c))^3}. 
$$

Lemma 16 now states that $c\in V_{k-2}\setminus V_{k-3}$ and we
can bound $\Sigma_{(2)}$,
distinguishing $k=0,1$ from $k\ge 2$  as:
$$
\Sigma_{(2)} \le 
\sum_{t=\tau}^{n} 
\left(
\sum_{\scriptsize
\begin{array}{c}
0\le k\le 1,\cr h=n-t-k
\end{array}
}
8\cdot \#V_k\cdot\frac{8^\beta}{t^{3\beta}} \right. 
\times
\left( \sum_{\Delta \in {\rm Til}_{h-1} } \left(27{\rm mes
}\Delta\right)^{\beta} \right)
+
$$
$$
\left.
+
\sum_{\scriptsize
\begin{array}{c}
k\geq 2,h\ge 0: \cr k+h+t=n
\end{array}
}
 \left( 
\sum_{c\in V_{k-2}\setminus V_{k-3}}
\frac{8}{(q(c))^{3 \beta}}\right)
\times \frac{8^\beta}{t^{3\beta}} 
\times
\left( \sum_{\Delta \in {\rm Til}_{h-1} } \left(27{\rm mes
}\Delta\right)^{\beta} \right)
\right)
\le 
$$
$$\le
\frac{216^\beta\cdot 8}{\tau^{3\beta-1}} \, \, \times \, \, 
\left(9+\sum_{k=2}^\infty
\left( \sum_{c\in V_{k-2}\setminus V_{k-3}}
\frac{1}{(q(c))^{3 \beta}}\right) \, \,  
\times \, \, 
\sum_{h=0}^\infty \left( \sum_{\Delta \in {\rm Til}_{h} }
\left({\rm mes }\Delta\right)^{\beta} 
\right)\right). 
$$
But
$$
\sum_{h=0}^\infty \left( \sum_{\Delta \in {\rm  Til}_{h} }  
\left({\rm mes }\Delta\right)^{\beta} \right) =
\sum_{h=0}^\infty 
\sigma_{h, \beta } (\hbox{\got B}) \le
\frac{32}{3}2^\beta\zeta(2\beta)\zeta(3\beta-2)  
$$
by Lemma 14, and
$$
9+\sum_{k=2}^\infty  \sum_{c\in V_{k-2}\setminus V_{k-3}}
\frac{1}{(q(c))^{3\beta}}
\le
9\ +\ 8
\sum_{q=1}^\infty
\sum
_{\scriptsize
\begin{array}{c}
1\le a,b,\le q \\
g.c.d.(a,b,q)=1
\end{array}
}
\frac{1}{q^{3\beta}}
\le 
\left(9+\frac{32}{3}\right)\zeta(3\beta-2),  
$$
applying the upper bound (\ref{zeta}). Now the inequality of Lemma 17
follows.\hfill$\Box$

{\bf Lemma 18.}\, \, \, 
{For all $\beta > 1$,
$$\Sigma_{(3)}=\frac{L(\text{\got B},3\beta)}{(n^2/2)^\beta}
+ O\left(\Sigma_{(1)}\right)+ O\left(\frac{1}{n^{2\beta} w^{3(\beta - 
1)}}+\frac{w}{n^{2\beta + 1}}\right)
$$
}

\proof  
W.l.o.g., let $n-w$ be even (otherwise use $w'=w+1$, covering even
more cases). Let 
$$
\Sigma_{(3)}' = \sum_{\scriptsize
\begin{array}{c}
\Delta \in {\rm Til}_n:\cr \co(\Delta)
\, \, \text{such}\, \, \text{that }\, \, \co_{w+1}=\dots=\co_n=0
\end{array}}
\left({\rm mes }\Delta\right)^\beta.
$$

Apparently,
$$\Sigma_{(3)}=\Sigma_{(3)}'+O(\Sigma_{(1)}).$$
Now, every triangle $\Delta\in\text{Til}_n$ with code
$\co(\Delta)=\co_1\dots \co_w0^{n-w}$ is a subset of a unique triangle
$\Delta'\in \text{Til}_w$ with code $\co(\Delta') = \co_1\dots \co_w$. 

Let $(a,b,c)$ and $(a',b',c')$ be the vertices of $\Delta$ and
$\Delta'$, respectively. Since $\Delta$ is obtained by $k = n-w$
(which is even) operations ``0'', with Lemma 2 and Lemma 13(iii), we
get
$$\text{mes}\Delta'=\frac{1}{2
\(q(a)+\frac{k}{2}q(c)\) \(q(b)+\frac{k}{2}q(c)\)q(c)}$$
and thus
$$
\sum_{c\in V_w}
\frac{\operatorname{deg}(c)}{\(\frac{n^2}{2}q(c)^3\)^\beta}
 \le 
\Sigma_{(3)}'
\le 
\sum_{c\in V_w}
\frac{\deg(c)}{\(\frac{(n-w+1)^2}{2}q(c)^3\)^\beta}.$$
Since
$$
\sum_{q=1}^{w-1} \frac{\sum_{l} lG_l(q)}{q^{3\beta}} \le 
\sum_{m=0}^{w-1}
 \sum_{ c \in V_{m}\setminus V_{m-1}}
 \frac{\deg(c)}{(q(c))^{3\beta}} = \sum_{ c \in V_{w}}
 \frac{\deg(c)}{(q(c))^{3\beta}} =L(\hbox{\got B}, 3\beta) +O\left(
 w^{-3(\beta - 1)} 
\right),
$$
as in (\ref{SumDeg1}), (\ref{SumDeg2}), it follows that
$$
\Sigma_{(3)}' = 
 \frac{L(\hbox{\got B}, 3\beta)}{(n^2/2)^\beta}  
\( 1+O\( \frac{w}{n}\)\), 
$$
and the lemma is proved.\hfill $\Box$

{\bf Theorem 3.}\, \, \, \, {\it For $\beta > 1$ the following
  asymptotic formula is valid 
$$
\sigma_{n, \beta } (\hbox{\got B}) = 
\frac{2^\beta\cdot L(\hbox{\got B},3\beta)}{n^{2\beta} }\left(1+ O\left( \frac{(\log 
n)^{1-\frac{1}{3\beta}}}{n^{\frac{\beta -1}{3\beta}}}\right)\right) .
$$
}

\proof
Assembling the results of Lemmata 15, 17, and 18, 
$$
\sigma_{n, \beta } (\hbox{\got B}) 
= 
\frac{L(\hbox{\got B}, 3\beta)}{(n^2/2)^{\beta} }
+ O\left( 
\frac{1}{n^{{\log 2\, \,  \gamma(\beta - 1)}}} 
+ \frac{1}{n^{2\beta} w^{3(\beta - 1)}}
+\frac{w}{n^{2\beta + 1}} 
+ \left( \frac{\gamma \log n}{w}\right)^{3\beta - 1} 
\right)=
 $$
 $$
= \frac{L(\hbox{\got B}, 3\beta)}{(n^2/2)^{\beta} }
\left(1+ O\left( 
\frac{(\log n)^{1-\frac{1}{3\beta}}}{n^{\frac{\beta -1}{3\beta}}}
\right)\right),
 $$
we obtain the
result.\hfill$\Box$ 


{\sc
Department of Theory of Numbers,\ 
Fac.~Mathematics and Mechanics,\ 
Moscow State University,\ 
Leninskie Gory,\ 
119992, Moscow,\ 
RUSSIA\\ 
{\tt moshchevitin@rambler.ru}
}
\\\\
{\sc
Instituto de Matem\'aticas,\ 
Facultad de Ciencias,\ 
Universidad Austral de Chile,\ 
Casilla 567,\ 
Valdivia,\  
CHILE\\  
{\tt vielhaber@gmail.com}
}

\begin{thebibliography}{77}
\bibitem{Stern}
M.~Stern, {\it \"Uber eine zahlentheoretische Funktion}, Crelles
Journal f\"ur die reine und angewandte Mathematik  55 (1858), 193--220. 
\bibitem{BR}
Brocot A., {\it Calcul des rouages par approximation, 
nouvelle m\'ethode} Revue Chronom\'etrique 6 (1860), 186--194.

\bibitem{LU} Lucas E., {\it Th\'eorie  des Nombres},
  Gauthiers--Villars, Paris,  Vol. 1 (1891), 467--475, 508--510.
\bibitem{MZ}   Moshchevitin N., Zhigljavsky A., 
{\it Entropies of the partitions of the unit interval
generated by the Farey tree},  Acta Arithmetica 115.1 (2004), 47--58.

\bibitem{DU}
Dushistova A., {\it On the partition of the unit interval generated by
  Brocot sequences}, to appear in: Sbornik Mathematics  198 No.~7
(2007) (Russian),  available at arXiv:math.NT/0512598v2 2 Jun 2006



\bibitem{SCHW} Schweiger F., {\it Multidimensional Continued Fraction
 Algorithms},  Oxford University Press, 2000.


\bibitem{BALADI} Baladi V., Nogueira A., {\it Lyapunov exponents for
  non-classical multidimensional continued fraction algorithms,} 
Nonlinearity  9 (1996), 1529--1546.

\bibitem{WOO}
Beaver O.V., Garrity T., {\it  A two-dimensional Minkowski ?(x)
  function},  J. Numb. Th.  107 (2004), 105--134. Available at
arXiv:math.Nt/0210480v2 8 Dec 2002


\bibitem{BRENTJES} Brentjes A., {\it Multidimensional continued fraction
algorithms},  Math. Centre Tract. 145 (1981).


\bibitem{HU}
Hurwitz A., {\it \"Uber die angen\"aherte Darstellung der Zahlen durch
rationale Br\"uche}, Math. Ann.  44 (1894),  417--436.

\bibitem{GRA}
Grabiner D.J., {\it Farey nets and multidimensional continued fractions},
Monatsh.~Math. 114 (1992), 35--60.


\bibitem{MON}
M\"onkemeyer R., {\it \"Uber Fareynetze in $n$ Dimensionen}, 
Math. Nachr. 11 (1963), 321--344.

\bibitem{PUA}
Poincar\'e H., {\it  Sur une g\'en\'eralisation des fractions
  continues},
C.R.~Acad.Sci.Paris  Ser. A 99 (1884),   1014--1016  = 
Oeuvres V, pp. 185--187.

\bibitem{NOUG}
Nogueira A., {\it
 The three-dimensional Poincar\'e continued fraction algorithm},
 Israel J. Math.  90 (1995),  373--401.

\end{thebibliography}
\end{document}